%% file: Royset_setconv.tex
\documentclass[titlepage,final,11pt]{article}
\input mac.tex

\begin{document}


\begin{center}
\begin{large}
{\bf Set-Convergence and Its Application: A Tutorial}
\smallskip
\end{large}
\vglue 0.3truecm
\begin{tabular}{c}
  \begin{large} {\sl Johannes O. Royset 
                                  } \end{large} \\
  Operations Research Department\\
  Naval Postgraduate School, Monterey, California\\
\end{tabular}

\vskip 0.1truecm

\end{center}

\vskip 0.3truecm

\noindent {\bf Abstract}. \quad Optimization problems, generalized equations, and the multitude of other variational problems invariably lead to the analysis of sets and set-valued mappings as well as their approximations. We review the central concept of set-convergence and explain its role in defining a notion of proximity between sets, especially for epigraphs of functions and graphs of set-valued mappings. The development leads to an approximation theory for optimization problems and generalized equations with profound consequences for the construction of algorithms. We also introduce the role of set-convergence in variational geometry and subdifferentiability with applications to optimality conditions. Examples  illustrate the importance of set-convergence in stability analysis, error analysis, construction of algorithms, statistical estimation, and probability theory.

\vskip 0.2truecm

\halign{&\vtop{\parindent=0pt
   \hangindent2.5em\strut#\strut}\cr
{\bf Keywords}:  set-convergence,  epi-convergence, graphical convergence, weak convergence, stability,\newline approximation theory, variational geometry, subdifferentiability, truncated Hausdorff distance.\\

{\bf Date}:\quad \ \today \cr}

\baselineskip=15pt

\section{Introduction}

In the study of optimization problems, we quickly encounter a plethora of sets that might cause concerns for the classically trained: Sets of solutions are rarely singletons, local properties aren't characterized by a single gradient vector but sets of subgradients, and the accuracy of function approximations is determined by the distance between certain sets. More general problems such as variational inequalities and generalized equations are even defined in terms of set-valued mappings. The mathematics of sets is therefore central to most analysis of optimization and variational problems beyond the elementary cases. In particular, a notion of convergence of sets becomes necessary, for example, to clarify what it means for a sequence of solution sets or a sequence of subgradient sets to converge as well as for many other concepts.

In 1902, Painlev\'{e} defined set-convergence in the sense we have it today, with Hausdorff and Kuratowski supporting its dissemination by including the topic in their books in 1927 and 1933, respectively\footnote{Detailed historical remarks can be found at the end of Chapter 4 in \cite{VaAn}.}. A set can be viewed as a point in a space of sets on which a topology and other mathematical structure can be defined. For the majority of the 20th century, set-convergence was studied from this point of view culminating with Beer's monograph \cite{Beer.93}. The growing number of variational problems in statistics, medicine, economics, business, engineering, and the sciences now drives a renewed interest in set-convergence from a more applied angle. In this article, we summarize key properties of set-convergence and highlight their importance when analyzing optimization problems and generalized equations as well as corresponding algorithms.

An illustrative example arises in the classical penalty method for constrained optimization. For $f_0,h:\reals^n\to \reals$, the problem of minimizing $f_0(x)$ subject to $h(x) = 0$ can be stated concisely as\footnote{For any set $C$, $\iota_C(u) = 0$ if $u\in C$ and $\iota_C(u) = \infty$ otherwise.}
\[
\nnmin_{x\in \reals^n} f(x) = f_0(x) + \iota_{\{0\}}\big(h(x)\big).
\]
Since an infinite penalty isn't compatible with direct application of standard algorithms, the penalty method considers the alternative problems
\[
\nnmin_{x\in \reals^n} f^\nu(x) = f_0(x) + \theta^\nu\big(h(x)\big)^2, ~\mbox{ with } \theta^\nu>0, ~~ \nu\in\nats = \{1, 2, \dots\}.
\]
At least for smooth $f_0$ and $h$, standard unconstrained algorithms apply to the alternative problems but can we guarantee that their minimizers are near those of the actual problem? Intuitively, this relates to whether $f^\nu$ is near $f$ in some sense. A classical treatment is complicated by the fact that $f(x) - f^\nu(x)=\infty$ whenever $h(x)\neq 0$ regardless of the value of $\theta^\nu$. Thus, it appears that the difference between $f$ and $f^\nu$ remains large even as $\theta^\nu\to \infty$ and one might be led to believe that the alternative problems aren't suitable approximations of the actual problem.

For another example, consider the problem of minimizing $f_0(x) = -x$ subject to $g(x) = x^3 - x^2 - x + 1\leq 0$; the feasible set is $(-\infty, -1] \cup \{1\}$ and $x^\star=1$ is the unique minimizer. If $g$ is replaced by the seemingly accurate approximation $g^\nu = g + \nu^{-1}$ so that $\nsup_{x\in\reals} |g^\nu(x) - g(x)| = \nu^{-1}$, then one might expect minimizers to change just a little. However, we obtain a substantially different minimizer near $-1$ even for large $\nu$. In these two examples, an assessment based on differences in function values is misleading. The key insight is provided by the epigraphs of the alternative functions and their proximity to that of the actual function.

For a function $f:\reals^n\to \Reals = [-\infty,\infty]$, we recall that the {\it epigraph} $\epi f = \{(x,\alpha)\in\reals^n\times \reals~|~f(x) \leq \alpha\}$. Moreover, the {\it minimum value} of $f$ is $\inf f = \inf \{f(x)~|~x\in\reals^n\}$ and the {\it set of minimizers} is $\nargmin f = \{x\in \dom f~|~f(x) = \inf f\}$, where $\dom f = \{x\in\reals^n~|~f(x)<\infty\}$.

\drawing{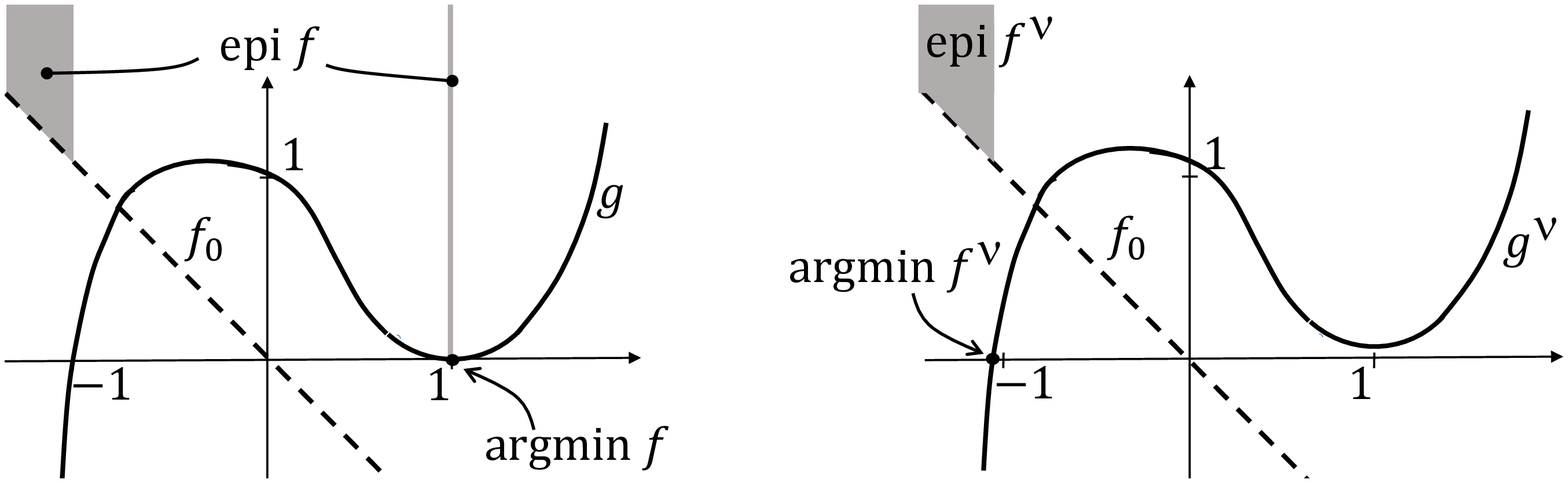}{4.2in} {Constraint function $g(x) =
x^3 - x^2 - x + 1$ is approached uniformly by $g^\nu = g +
\nu^{-1}$, but minimum value and minimizer of $f^\nu = f_0 + \iota_{(-\infty, 0]}(g^\nu(\cdot))$ aren't close to those of $f = f_0 + \iota_{(-\infty, 0]}(g(\cdot))$.}
{fig:unif}

In the second example above, we compare the actual problem of minimizing $f = f_0 + \iota_{(-\infty, 0]}(g(\cdot))$ with that of minimizing $f^\nu = f_0 + \iota_{(-\infty, 0]}(g^\nu(\cdot))$. Figure \ref{fig:unif} shows $\epi f$ and $\epi f^\nu$. The two epigraphs are rather different and they don't get close in any reasonable sense; $\epi f$ contains a vertical ray starting at $(1,0)$ that is distinct from any points in $\epi f^\nu$. The effect is that $\nargmin f^\nu$ and $\inf f^\nu$ deviate greatly from $\nargmin f$ and $\inf f$, respectively. Indeed, the proximity between epigraphs provides a robust assessment of whether the corresponding solutions are close as we see below.

To formalize set-convergence, we define the {\it point-to-set
distance}\index{distance!point-to-set} between $\bar x\in \reals^n$ and
$C\subset \reals^n$ as
\begin{equation}\label{eqn:distdef}
  \dist_2(\bar x,C) = \ninf_{x\in C} \|x-\bar x\|_2 ~\mbox{ when } C\neq \emptyset~ \mbox{ and } \dist_2(\bar x,\emptyset) = \infty.
\end{equation}

\begin{definition}\label{dSetLimit}{\rm (set-convergence\footnote{This notion is also referred to as Painlev\'{e}-Kuratowski convergence.}).}
For $\{C, C^\nu\subset \reals^n, ~\nu\in\nats\}$, we say that $C^\nu$ set-converges to $C$, written $C^\nu\sto C$ or $\nLim C^\nu = C$, when
\[
C \mbox{ is closed~ and~ } \dist_2(x,C^\nu)\to \dist_2(x,C) ~\mbox{ for all } x\in \reals^n.
\]
\end{definition}

\drawing{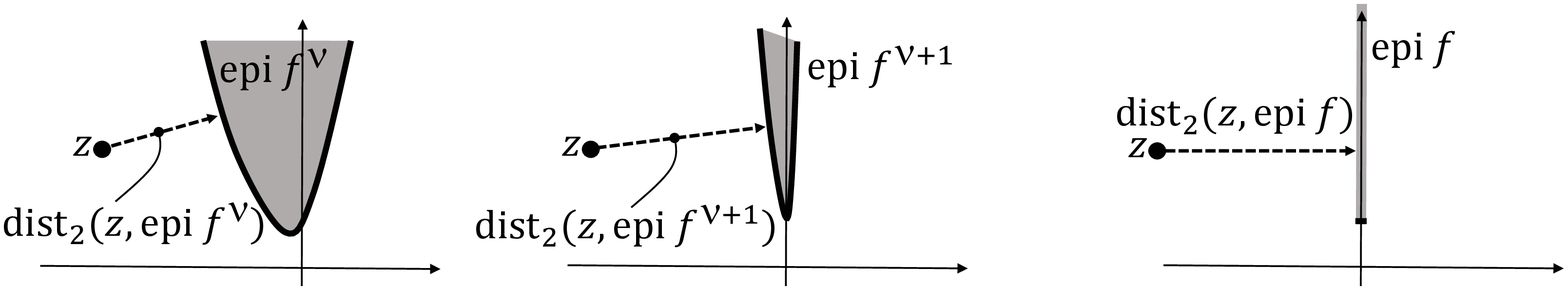}{5.5in} {The epigraphs of $f^\nu$ set-converge to $\epi f = \{0\} \times [1, \infty)$.}{fig:epiconvPenalty}

Figure \ref{fig:epiconvPenalty} visualizes the definition applied to $\epi f^\nu$ and $\epi f$ in the first example above when $f_0(x) = (x+1)^2$ and $h(x) = x$. Indeed, $\dist_2(z,\epi f^\nu)\to \dist_2(z,\epi f)$ for all $z\in \reals^{n+1}$. As expected from the penalty method, we see that $\nargmin f^\nu$, which is a bit to the left of the origin, tends to $\nargmin f = \{0\}$.

In general, for any $f^\nu:\reals^n\to \Reals$ and a proper\footnote{A function $f:\reals^n\to \Reals$ is {\it proper} if $f(x)>-\infty$ for all $x\in\reals^n$ and $f(x)<\infty$ for some $x\in\reals^n$.} function $f:\reals^n\to \Reals$:
\[
\epi f^\nu\sto \epi f \mbox{ and $\bar x$ is a cluster point of } \{x^\nu \in \nargmin f^\nu, \nu\in\nats\} ~~\Longrightarrow~~ \bar x\in \nargmin f.
\]

This fact, further elaborated on in \S\ref{sec:optim}, guides the construction of algorithms such as the penalty method and helps us reformulate problems that are fundamentally unstable. As an example, consider
\[
  \nnmin_{x\in\reals^n}  f_0(x)    \mbox{ subject to }  g_i(x) \leq 0, ~i=1, \dots,q,
\]
where the individual functions need to be replaced by approximations for computational reasons or because they are only partially known. Regardless, let $f_0^\nu, g_i^\nu, \dots, g_q^\nu$ be the approximations and, for $\alpha^\nu\geq 0$, suppose that $\sup_{x\in\reals^n} |f_0^\nu (x) - f_0(x)| \leq \alpha^\nu$ and $\sup_{x\in\reals^n} \max_{i=1,\dots, q} |g^\nu_i (x) - g_i(x)| \leq \alpha^\nu$, with all the functions being real-valued. As discussed in conjuncture with Figure \ref{fig:unif}, one can't expect that substituting in $f^\nu_0$ for $f_0$ and $g^\nu_i$ for $g_i$ produce a problem with minimizers near those of the actual one even if $\alpha^\nu$ is small. The goal becomes to construct alternative functions with epigraphs that set-converge to one corresponding to the actual problem as $\alpha^\nu\to 0$ because then minimizers will follow suit.

An alternative and indeed approximating problem emerges by softening the constraints and using a penalty $\theta^\nu>0$:
\[
\nnmin_{x\in\reals^n, y\in\reals^q} ~f^\nu_0(x) + \theta^\nu \nsum_{i=1}^q y_i  ~\mbox{ subject to }~  g^\nu_i(x) \leq y_i, ~~0 \leq y_i, ~~i=1, \dots,q.
\]
Although seemingly rather different than the actual problem, it's epigraphically close in the sense we seek. Specifically, the approximating problem can be written compactly as minimizing the function $f^\nu:\reals^n\times \reals^q\to \Reals$ given by
\[
f^\nu(x,y) = f^\nu_0(x)  + \nsum_{i=1}^q \phi^\nu(y_i) + \nsum_{i=1}^q \iota_{(-\infty, 0]}\big(g^\nu_i(x) - y_i\big), \mbox{ where } \phi^\nu(\beta) = \begin{cases}
\theta^\nu \beta & \mbox{  if } \beta\geq 0\\
\infty & \mbox{ otherwise.}
\end{cases}
\]
The actual problem is expressed in terms of $x$ only, but to allow comparison with the approximation we artificially introduce $y$ and state it as minimizing the function $f:\reals^n\times \reals^q\to \Reals$ defined by
\[
f(x,y) = f_0(x) + \nsum_{i=1}^q \iota_{\{0\}}(y_i) + \nsum_{i=1}^q \iota_{(-\infty, 0]}\big(g_i(x) - y_i\big).
\]
Trivially, $y$ must be the zero vector for $f$ to be finite so this is indeed an equivalent statement of the actual problem. The similarities between $f^\nu$ and $f$ are clear. In particular, $\epi \phi^\nu \sto \epi \iota_{\{0\}}$ when $\theta^\nu\to \infty$. Using the techniques of \S\ref{sec:optim}, we can show that $\epi f^\nu\sto \epi f$ under the assumptions that $f_0$ is continuous, $g_i$ is lower semicontinuous (lsc) for all $i$, $\theta^\nu\to \infty$, and $\theta^\nu \alpha^\nu\to 0$. Practitioners often soften constraints and the above analysis reveals that there are good reasons for this: the resulting problem becomes a valid approximation of an actual problem, which may not even be fully known. In contrast, the naive ``approximation'' obtained by simply substituting in $f_0^\nu$ and $g_i^\nu$ can be arbitrarily poor.

These introductory examples illustrate the role of set-convergence in analysis of optimization problems and the construction of numerical procedures for their solution. After some background on set-convergence in \S\ref{sec:found}, we return to that subject in \S\ref{sec:optim}. Set-convergence also enters in the development of robust notions of normal cones and subgradients that reach beyond the convex case as seen in \S\ref{sec:vargeom}. Continuity and approximation of set-valued mappings are intimately tied to set-convergence, with significant implications for the solution of generalized equations as discussed in \S\ref{sec:geneq}. Quantification of set-convergence enables us to obtain rates of convergence for algorithms and bound the solution error due to approximations in optimization problems and generalized equations; see \S\ref{sec:quant}. Applications of set-convergence in probability and statistics are recorded in \S\ref{sec:optim}. The article ends in \S\ref{sec:concl} with a summary of topics not covered including extensions to infinite dimensions.

\section{Inner and Outer Limits}\label{sec:found}

Set-convergence can be understood and fully characterized using the notions of outer and inner limits, which in fact was the approach of Painlev\'{e}. These are not only useful stepping stones towards set-convergence, but also furnish intermediate properties that might suffice in a particular application.

The {\it inner limit} of $\{C^\nu \subset\reals^n, \nu\in\nats\}$, denoted by $\nInnLim C^\nu$, is the collection of limit points to which sequences of points selected from the sets converge. Specifically,
\begin{equation*}
  \nInnLim C^\nu = \big\{x\in \reals^n~\big|~\exists x^\nu\in C^\nu\to x\big\}.
\end{equation*}
For example, the inner limit of the line segments $C^\nu = [1/\nu, 2/\nu]$ is $\{0\}$. A slight clarification regarding empty sets is in place. Suppose that $C^1 = \emptyset$, but $C^\nu = [1/\nu, 2/\nu]$ for $\nu \geq 2$. Still, $\nInnLim C^\nu = \{0\}$. Since convergence of a sequence is always determined by its tail and not a finite number of its elements, we aren't concerned about not being able to select a point from $C^1$. It suffices in the definition of inner limits to have points in $C^\nu$ converging to $x$ for all $\nu$ sufficiently large.

Figure \ref{fig:setconv2} furnishes two more examples: On the left, $\nInnLim C^\nu$ is the third quadrant. On the right, collections of vertical line segments $C^\nu$ eventually ``fill'' a rectangle.

\drawing{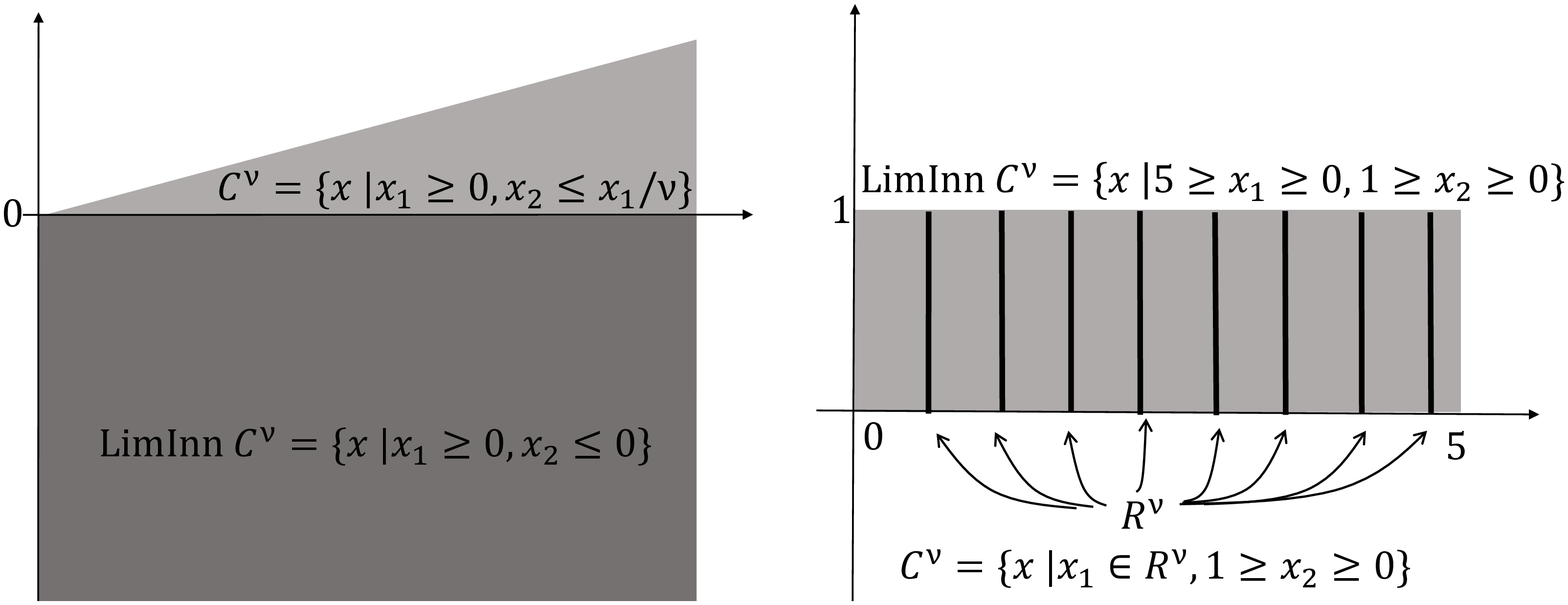}{4.8in}
   {Examples of inner limits; $R^\nu$ is the first $\nu$ numbers in an ordered list of the rational numbers in $[0,5]$.} {fig:setconv2}

If $C^\nu = \{0\}$ when
$\nu$ is odd and $C^\nu = [0, 1]$ when $\nu$ is even, then $\nInnLim
C^\nu=\{0\}$. The inner limit isn't $[0, 1]$ because we can't
construct $x^\nu \in C^\nu$ converging to $x\in (0,1]$.
This brings us to outer limits.

The {\it outer limit} of $\{C^\nu \subset\reals^n, \nu\in\nats\}$, denoted by $\nOutLim C^\nu$, is the collection of cluster points to which subsequences of points selected from the sets converge. To write this formally, we use the following notation for subsequences:
\[
\mbox{$\cN_\infty^\grill$ is the set of all infinite collections of increasing numbers from $\nats$}.\]
For example, $\{1, 3, 5, \dots\} \in \cN_\infty^\grill$ and $\{2, 4, 6, \dots\}\in \cN_\infty^\grill$. A subsequence of $\{x^\nu\in \reals^n, \nu\in\nats\}$ is then of the form $\{x^\nu\in \reals^n, \nu\in N\}$ for some $N\in \cN_\infty^\grill$, with $x^\nu \Nto x$ expressing that $x$ is a cluster point corresponding to the subsequence specified by the index set $N$. In this notation,
\begin{equation*}
  \nOutLim C^\nu = \big\{x\in \reals^n~\big|~ \exists N\in \cN_\infty^\grill \mbox{ and } x^\nu\in C^{\nu} \Nto x\big\}.
\end{equation*}
Hence, $x\in \nOutLim C^\nu$ if we can select an index set $N\in \cN_\infty^\grill$ and points $\{x^\nu\in C^\nu, \nu\in N\}$ converging to $x$. In the example $C^\nu = \{0\}$ when $\nu$ is odd and $C^\nu =
[0, 1]$ when $\nu$ is even, $\nOutLim C^\nu =
[0, 1]$ because for any $x \in [0, 1]$, we can rely on $N = \{2, 4, \dots\}$ and select $x^\nu = x$ for $\nu\in N$. The relation between the inner and the outer limits in this example points to a fact that motivates the terminology ``inner'' and ``outer:'' Generally, $\nInnLim C^\nu \subset \nOutLim C^\nu$. Although $\nLim C^\nu$ may not exist (as is clear from the ``odd-even'' example above), the inner limit and the outer limit always exist but could be the empty set. Moreover, they are closed sets and characterize set-convergence \cite[Section 4.B]{VaAn}:

\begin{proposition}\label{pSetConvergence}{\rm (inner and outer limits).}
For $\{C, C^\nu\subset \reals^n, \nu\in\nats\}$,
\[
\nLim C^\nu = C ~~\Longleftrightarrow~~ \nInnLim C^\nu = \nOutLim C^\nu = C.
\]
\end{proposition}

In applications, a sequence of sets may not set-converge but useful insight can come from the inner and/or outer limits. For example, given a closed feasible set  $C\subset\reals^n$ that is approximated by $\{C^\nu\subset\reals^n, \nu\in\nats\}$, suppose that an algorithm  generates $\{x^\nu\in C^\nu, \nu\in\nats\}$, which we hope are near a feasible point in $C$. Indeed, if $\nOutLim C^\nu\subset C$, then every cluster point of $\{x^\nu, \nu\in\nats\}$ is contained in $C$ and is thus feasible for the actual problem. This might suffice as justification for the approximating sets $C^\nu$. If $C^\nu\sto C$, then we have in addition that every $x\in C$ can be approached by points in $C^\nu$, i.e., there is $x^\nu\in C^\nu \to x$.

Practical verification of set-convergence is supported by a collection of rules that enable us to conclude that images, inverse images, intersections, unions, sums, products, and other operations on sets preserve set-convergence, at least in part, under various assumptions. We refer to \cite[Section 4.G]{VaAn} for a detailed treatment and here concentrate on the central cases of intersections and unions.

A feasible set in an optimization problem is often specified as the intersection and/or union of a collection of other sets, each representing some requirement, which may need approximation. The question then becomes: how do approximations of the ``component'' sets affect the feasible set?

\begin{proposition}{\rm (union and intersection under set-convergence).}\label{pUnionInterSetConvergence}
 For sets $\{C^\nu, D^\nu \subset \reals^n, \nu\in\nats\}$,
\begin{align*}
\nLim \big( C^\nu \cup D^\nu \big) &\sto \nLim C^\nu \cup \nLim D^\nu\\
\nOutLim \big(C^\nu \cap D^\nu\big) &\subset \nOutLim C^\nu \cap \nOutLim D^\nu.
\end{align*}
If all $C^\nu, D^\nu$ are convex and\footnote{We denote by $\nt C$ the interior of a set $C\subset \reals^n$.} $\nt(\nLim C^\nu) \cap \nLim D^\nu \neq \emptyset$, then $\nLim ( C^\nu \cap D^\nu ) \sto \nLim C^\nu \cap \nLim D^\nu$.
\end{proposition}
\state Proof. The first two statements follow readily from the definitions of inner limits and outer limits. For the convex case, consult \cite[Theorem 4.32]{VaAn}.\eop

The situation is less than ideal for intersections: $\nOutLim(C^\nu \cap D^\nu)$ can be strictly contained in $\nLim C^\nu \cap \nLim D^\nu$ even for convex sets. For example, $C^\nu = [-1, -1/\nu]$ and $D^\nu = [1/\nu, 1]$ has $\nOutLim (C^\nu \cap D^\nu) =\emptyset$, while $\nOutLim C^\nu \cap \nOutLim D^\nu = \{0\}$.  This highlights the challenge in constrained optimization already indicated in the Introduction: A tiny change of the components of a problem may have large effects on the feasible set. Set-convergence helps us to identify trouble. If a feasible set $C \cap D$ isn't stable under perturbations of $C$ and $D$ in the sense that approximations $C^\nu \cap D^\nu$ fail to set-converge to $C\cap D$ despite $C^\nu\sto C$ and $D^\nu\sto D$, then optimization over $C\cap D$ is somehow ill-posed. Any small misspecification of $C$ or $D$ can cause misleading solutions. In the convex case, we achieve such stability if $C$ and $D$ overlap ``sufficiently'' as stated in the proposition.

Unions of sets are well-behaved in this sense, which has the following ramification for min-functions.

\begin{example}{\rm (approximation of min-functions).}
  For $f_i, f_i^\nu:\reals^n\to \Reals$, $i=1, \dots, m$, let
  \[
  f(x) = \min_{i=1, \dots, m} f_i(x) ~\mbox{ and }~ f^\nu(x) = \min_{i=1, \dots, m} f^\nu_i(x), ~~~\forall x\in\reals^n.
  \]
  Then, $\epi f^\nu \sto \epi f$ whenever $\epi f_i^\nu\sto \epi f_i$ for all $i$.
\end{example}
\state Detail. Since $\epi f = \cup_{i=1}^m \epi f_i$ and likewise for $\epi f^\nu$,  Proposition \ref{pUnionInterSetConvergence} establishes the conclusion.\eop

\section{Minimization and Epi-Convergence}\label{sec:optim}

Set-convergence of epigraphs goes back to the pioneering work by Wijsman \cite{Wijsman.64,Wijsman.66} and Mosco \cite{Mosco.69} on convex functions. With an increasing focus on the nonconvex case in the late 1970s, Wets coined the term epi-convergence\footnote{Independently, in their work on calculus of variation, Di Giorgi and co-authors came to the same notion of set-convergence of epigraphs and called it $\Gamma$-convergence; see \cite{DegioriFranzoni.75}.} \cite{Wets.80}: For functions $f, f^\nu:\reals^n\to \Reals$, $f^\nu$ {\it epi-converges} to $f$, written $f^\nu \eto f$, when the corresponding epigraphs set-converge, i.e.,
\[
f^\nu \eto f ~~~\Longleftrightarrow~~ \epi f^\nu \sto \epi f.
\]

Since inner limits are contained in the corresponding outer limits, $\epi f^\nu \sto \epi f$ takes place if and only if $\nOutLim (\epi f^\nu) \subset \epi f \subset \nInnLim (\epi f^\nu)$. The two inclusions translate into conditions (a) and (b) below; see \cite[Proposition 7.2]{VaAn}.

\begin{theorem}
\label{tEpiCnvr}{\rm (characterization of epi-convergence)}.
For $f,f^\nu:\reals^n\to \Reals$, $f^\nu\eto f$ if and only if
\begin{enumerate}[{\rm (a)}]
\vspace{-.4\baselineskip}\item $\forall x^\nu \to x$, $\nliminf f^\nu(x^\nu) \geq f(x)$,
\vspace{-.4\baselineskip}\item $\forall x, \,\exists\, x^\nu
    \to x$ such that $\nlimsup f^\nu(x^\nu) \leq f(x)$.
\end{enumerate}%
\end{theorem}

For (a), $f^\nu$ must be ``high enough'' regardless of how we approach $x$. For (b), there must be a ``path'' to $x$ along which $f^\nu$ is ``low enough''. Figure \ref{fig:epiconv1} illustrates this with $f(x) = -1$ for $x\leq 0$ and $f(x) = 1$ otherwise, and $f^\nu(x) = -1$ for $x \leq -1/\nu$, $f^\nu(x) = -\sqrt{-\nu x}$ for $x\in (-1/\nu, 0]$, $f^\nu(x) = \sqrt{\nu x}$ for $x\in (0,1/\nu]$, and $f^\nu(x) = 1$ otherwise. To verify epi-convergence, let's first consider (a). For any $x^\nu\to x\leq 0$, $f^\nu(x^\nu) \geq -1 = f(x)$ for all $\nu$. For any $x^\nu \to x>0$, $f^\nu(x^\nu) = 1 = f(x)$ for sufficiently large $\nu$. Thus, (a) holds. Second, let's consider (b). For $x\neq 0$, we can take $x^\nu = x$ because $f^\nu(x^\nu) = f^\nu(x) \to f(x)$. It only remains to show that for some $x^\nu \to 0$, $\nlimsup f^\nu(x^\nu) \leq f(0) = -1$. Here, we need to be careful: $x^\nu$ above 0 or too close to 0 has $f^\nu(x^\nu)$ much above $-1$. However, $x^\nu = -1/\nu$ results in $f^\nu(x^\nu) = -1$ for all $\nu$ and we still have $x^\nu\to 0$ so (b) holds.

\drawing{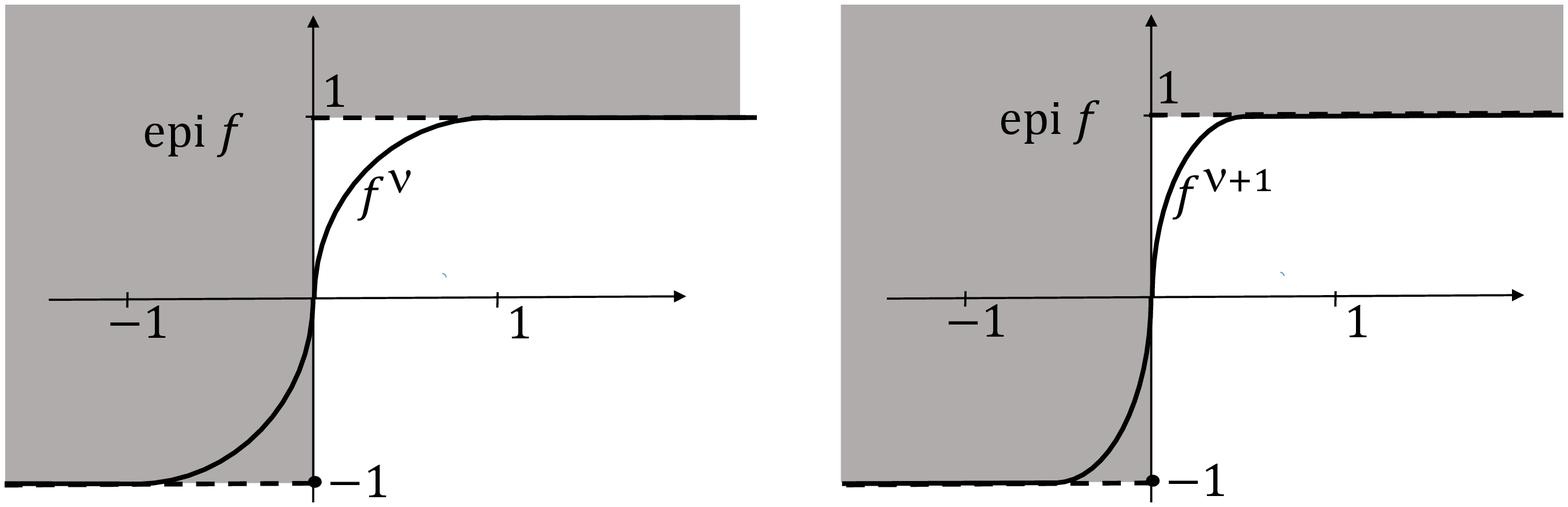}{4.2in} {Functions $f^\nu$ (solid line) epi-converge to $f$ (dashed line).}{fig:epiconv1}

\begin{example}\label{xWeakConv}{\rm (weak convergence of distribution functions).}
In probability theory, a sequence of distribution functions $F^\nu:\reals^n\to [0,1]$ converge weakly to a distribution function $F:\reals^n\to [0,1]$ when
\[
F^\nu(x) \to F(x) \mbox{ for all } x\in\reals^n \mbox{ at which } F \mbox{ is continuous.}
\]
This is an example of set-convergence. Specifically, by \cite[Theorem 3.2]{RoysetWets.16b},
\[
F^\nu(x) \mbox{ converges weakly to } F(x) ~~\Longleftrightarrow~~ -F^\nu \eto -F ~~\Longleftrightarrow~~ \hypo F^\nu \sto \hypo F,
\]
where for any $f:\reals^n\to \Reals$, $\hypo f = \big\{(x,\alpha)\in\reals^n\times \reals~\big|~f(x) \geq \alpha\big\}$ is the hypograph\footnote{In general, set-convergence of hypographs is referred to as {\it hypo-convergence} with analogous properties to those of epi-convergence but with its primary application area being maximization instead of minimization.} of $f$.
\end{example}
\state Detail. It's the monotonicity and upper-semicontinuity of distribution functions that allow us to make the equivalence between pointwise convergence and epi-convergence, which doesn't hold in general. If Figure \ref{fig:epiconv1} is modified by setting $f(0) = 1$ and $f$ is replaced by $\half f + \half$ then the resulting function, say $F$, is a distribution function with a discontinuity at $x=0$. Likewise, if we replace $f^\nu$ by $\half f^\nu + \half$, then the resulting function $F^\nu$ is a distribution function. For $x\neq 0$, $F^\nu(x)\to F(x)$ so $F^\nu$ converges weakly to $F$. Meanwhile, using a similar argument as in the discussion above, we obtain that $-F^\nu\eto -F$.\eop

\begin{example}\label{xKWmaximumlikelihood}{\rm (density estimation using location mixtures).} Given a sample $\xi^1, \dots, \xi^n\in\reals^m$, the Kiefer-Wolfowitz maximum likelihood estimate of a probability density function on $\reals^m$, used in clustering and denoising, is a minimizer of
\[
\nnmin_{p\in \cP} -\frac{1}{n}\nsum_{j=1}^n \ln p(\xi^j),
\]
where $\cP$ is the family of location mixtures of the $m$-variate standard normal density $\phi = (2\pi)^{-m/2} \exp(-\half\|\cdot\|_2^2)$; see \cite{SahaGuntuboyina.20}. The problem is equivalent to minimizing $-\frac{1}{n}\sum_{j=1}^n \ln x_j$ subject to\footnote{For $C\subset\reals^n$, $\con C$ is its convex hull.}
\[
x\in C = \con\Big\{\big(\phi(\xi^1-z), \dots, \phi(\xi^n-z)\big)~\Big|~ z\in\reals^m\Big\}.
\]
Since $C$ is the convex hull of an infinite number points, it needs to be approximated.
\end{example}
\state Detail. Given $\bar z^1, \dots, \bar z^\nu\in \reals^m$, we can approximate $C$ by 
\[
C^\nu = \con\Big\{\big(\phi(\xi^1-\bar z^k), \dots, \phi(\xi^n-\bar z^k)\big)~\Big|~ k=1, \dots, \nu\Big\}.
\]
The approximating problem with $C$ replaced by $C^\nu$ is tractable because of
the equivalent convex problem
\[
\nnmin_{\mu\in\reals^\nu} -\frac{1}{n}\nsum_{j=1}^n \ln \left(\nsum_{k=1}^\nu \mu_k \phi(\xi^j-\bar z^k) \right) \mbox{ s.t. } \nsum_{k=1}^\nu \mu_k = 1, ~\mu_k\geq 0~ \forall k=1,\dots, \nu.
\]

To justify the approximation, let $f(x) = -\frac{1}{n}\sum_{j=1}^n \ln x_j + \iota_C(x)$ and $f^\nu(x) = -\frac{1}{n}\sum_{j=1}^n \ln x_j + \iota_{C^\nu}(x)$, which represent the actual and approximating problems, respectively. (Here, $-\ln \alpha = \infty$ for $\alpha \leq 0$.) We show that $f^\nu\eto f$ when $\{\bar z^k,
k\in\nats\}$ has the property that for any $\bar z\in \reals^m$, there exists $\bar z^k\to \bar z$.

First, consider  Theorem \ref{tEpiCnvr}(a) and let $x^\nu \to x$. Since $C^\nu\subset C$ and $f$ is lsc, $\nliminf f^\nu(x^\nu)\geq \nliminf f(x^\nu)\geq f(x)$. Second, consider Theorem \ref{tEpiCnvr}(b) and let
$x\in C$. We need to construct $x^\nu\in C^\nu\to x$. By
Caratheodory's Theorem, there exists $\alpha_0, \dots.,
\alpha_n\geq 0$ and $z^0, \dots, z^n\in\reals^m$ such that $\sum_{j=0}^n
\alpha_j = 1$ and $x = \sum_{j=0}^n \alpha_j (\phi(\xi^1-z^j), \dots.,
\phi(\xi^n-z^j))$. Set $x^\nu = \nsum_{j=0}^n
\alpha_j (\phi(\xi^1-u_j^\nu), \dots., \phi(\xi^n-u_j^\nu))$,  with $u_j^\nu\in\nargmin_{k=1, \dots, \nu} \|z^j - \bar z^k\|_2$. Certainly, $x^\nu\in C^\nu$. Since $\bar z^1, \bar z^2, \dots$ eventually
``fill'' $\reals^m$ by assumption, $u_j^\nu\to z^j$ as $\nu\to \infty$ for all $j$. By the continuity of $\phi$, we have that $x^\nu\to x$ and also
$f^\nu(x^\nu)\to f(x)$ by the continuity of the log-function. We have
established the conditions of Theorem  \ref{tEpiCnvr} and $f^\nu\eto f$.\eop

For other applications of set-convergence in probability and statistics; see \cite{RockafellarRoyset.14,Aswani.19} and \S\ref{sec:concl}.

Epi-convergence of composite function can often be established by checking properties of its various components. This might offer an easier path than working directly with Theorem \ref{tEpiCnvr}.

\begin{proposition}
\label{pEpiCalculusSum}{\rm (epi-convergence of sums)}.
For functions $f,f^\nu,g,g^\nu:$ $\reals^n\to (-\infty, \infty]$, suppose that $f^\nu\eto f$ and $g^\nu\eto g$. Then,
$f^\nu + g^\nu \eto f + g$ under either one of the following two conditions:
\begin{enumerate}[{\rm (a)}]
\vspace{-.4\baselineskip}\item $-g^\nu \eto -g$.

\vspace{-.4\baselineskip}\item $g^\nu(x) \to g(x)$ and $f^\nu(x) \to f(x)$ for all $x\in\reals^n$.
\end{enumerate}
The pointwise convergence of $g^\nu$ ($f^\nu$) in (b) holds automatically if all $g^\nu$ ($f^\nu$) are convex and real-valued.
\end{proposition}
\state Proof. For sequences $\{\alpha^\nu, \beta^\nu\in\Reals, \nu\in\nats\}$, $\nliminf (\alpha^\nu + \beta^\nu) \geq \nliminf \alpha^\nu + \nliminf \beta^\nu$ provided that the right-hand side is not $\infty-\infty$ or $-\infty + \infty$, and likewise for limsup but then with $\leq$. If $x^\nu\to x$, then
\[
\nliminf \big(f^\nu(x^\nu) + g^\nu(x^\nu)\big) \geq \nliminf f^\nu(x^\nu) + \nliminf g^\nu(x^\nu) \geq f(x) + g(x)
\]
and Theorem \ref{tEpiCnvr}(a) holds without either condition (a) or (b). Given $x\in\reals^n$, there exists $x^\nu\to x$ such that $\nlimsup f^\nu(x^\nu) \leq f(x)$. We note that $\nlimsup f^\nu(x^\nu)>-\infty$ because otherwise $\nliminf f^\nu(x^\nu) = -\infty$ and $f(x) = -\infty$, but this violates the fact that $f>-\infty$. Under (a), $g^\nu(x^\nu) \to g(x)$ and
\[
\nlimsup \big(f^\nu(x^\nu) + g^\nu(x^\nu)\big) \leq \nlimsup f^\nu(x^\nu) + \nlimsup g^\nu(x^\nu) \leq f(x) + g(x)
\]
and Theorem  \ref{tEpiCnvr}(b) holds. Under condition (b), set $x^\nu = x$ so that
\[
\nlimsup \big(f^\nu(x^\nu) + g^\nu(x^\nu)\big) \leq \nlimsup f^\nu(x) + \nlimsup g^\nu(x) = f(x) + g(x)
\]
and Theorem \ref{tEpiCnvr}(b) holds again.

When all the functions $g^\nu$ are convex and real-valued, then $g^\nu\eto g$ if and only if $g^\nu(x)\to g(x)$ for all $x\in\reals^n$ by \cite[Theorem 7.17]{VaAn}. Thus, (b) holds automatically under the additional assumption.\eop

We note that having $g^\nu\eto g$ and $-g^\nu \eto -g$ is equivalent to having $g^\nu(x^\nu) \to g(x)$ for all $x^\nu\to x$. Thus, condition (a) of the proposition holds, for example, when $g^\nu = g$ and $g$ is continuous. It is now possible to verify the claim in the Introduction regarding epi-convergence for approximations obtained by constraint softening using Theorem \ref{tEpiCnvr} and Proposition \ref{pEpiCalculusSum}.

\begin{proposition}
\label{pEpiCalculusComposite}{\rm (epi-convergence of composite functions)}.
For $h,h^\nu:\reals\to\Reals$ and $f,f^\nu:\reals^n\to \Reals$, the following hold:
\begin{enumerate}[{\rm (a)}]

\vspace{-.4\baselineskip}\item If $h^\nu\eto h$, $h^\nu(\alpha) \to h(\alpha)$ for all $\alpha\in\reals$, and $f$ is real-valued and continuous, then $h^\nu \circ f \eto h \circ f$.

\vspace{-.4\baselineskip}\item If $f^\nu\eto f$ and $h$ is continuous, nondecreasing and extended to  $\Reals$ with the conventions $h(-\infty) = \inf_{\alpha\in\reals} h(\alpha)$ and $h(\infty) = \sup_{\alpha\in\reals} h(\alpha)$, then $h \circ f^\nu \eto h\circ f$.

\end{enumerate}
\end{proposition}
\state Proof. For (a), let $x^\nu\to x$, which ensures that $f(x^\nu)\to f(x)$. By Theorem  \ref{tEpiCnvr}(a), $\nliminf h^\nu(f(x^\nu)) \geq h(f(x))$. Given $x\in\reals^n$, set $x^\nu = x$. Then, $\nlimsup h^\nu(f(x^\nu)) = \nlimsup h^\nu(f(x)) = h(f(x))$. We have confirmed both conditions in Theorem  \ref{tEpiCnvr}.

For (b), let $x^\nu\to x$, which implies that $\nliminf f^\nu(x^\nu) \geq f(x)$ by  Theorem \ref{tEpiCnvr}(a).
Fix $\nu$ and let $\epsilon>0$. Suppose that $\alpha^\nu = \inf_{\mu\geq \nu} h(f^\mu(x^\mu))\in\reals$. Then, there exists $\bar \mu\geq \nu$ such that $\alpha^\nu \geq h(f^{\bar \mu}(x^{\bar \mu})) - \epsilon \geq h(\inf_{\mu\geq \nu} f^\mu(x^\mu)) - \epsilon$, the last inequality holds because $h$ is nondecreasing. Since $\epsilon$ is arbitrary, $\alpha^\nu\geq h(\inf_{\mu\geq \nu} f^\mu(x^\mu))$. A similar argument leads to the same inequality if $\alpha^\nu = -\infty$ and, trivially, also when $\alpha^\nu = \infty$.
Since the inequality holds for all $\nu$, it follows by the continuity of $h$ that
\begin{align*}
\nliminf h(f^\nu(x^\nu)) & = \lim_{\nu\to\infty} \big(\ninf_{\mu\geq \nu} h(f^\mu(x^\mu))\big) \geq \lim_{\nu\to\infty} h\big(\ninf_{\mu\geq \nu} f^\mu(x^\mu)\big)\\
& = h\big( \lim_{\nu\to\infty} ( \ninf_{\mu\geq \nu} f^\mu(x^\mu))\big) = h\big(\nliminf f^\nu(x^\nu)\big) \geq h(f(x)).
\end{align*}
Given $x\in\reals^n$, there exists $x^\nu\to x$ such that $f^\nu(x^\nu)\to f(x)$ by  Theorem \ref{tEpiCnvr}(b). Since $h$ is continuous, this implies  $h(f^\nu(x^\nu)) \to h(f(x))$ and in view Theorem \ref{tEpiCnvr} we have established (b).\eop

As indicated in the Introduction, $f^\nu \eto f$ ensures that cluster points of minimizers of $f^\nu$ are minimizers of $f$ but  the situation for minimum values is more delicate. If $f^\nu(x) = \max\{-1, x/\nu\}$ and $f(x) = 0$ for all $x$, then $f^\nu\eto f$ but $\inf f^\nu = -1 < \inf f= 0$. The trouble is that $\nLim (\nargmin f^\nu)=\emptyset$. This pathological situation is eliminated under {\it tightness.}

\begin{definition}{\rm (tightness)}.
The functions $\{f^\nu:\reals^n\to\Reals, \nu\in\nats\}$ are {\it tight} if for all $\epsilon >0$, there is a compact $B_\epsilon \subset \reals^n$ and an index $\nu_\epsilon$ such that
\[
    \ninf_{x\in B_\epsilon} f^\nu(x) \leq \ninf f^\nu +\epsilon \mbox{ for all } \nu\geq \nu_\epsilon.
\]
The functions {\em epi-converge tightly} if in addition to being tight they also epi-converge to some function.
\end{definition}

Tightness holds, for example, if $\{\dom f^\nu, \nu\in\nats\}$ are contained in
a bounded set or if there is a bounded set
$B\subset \reals^n$ such that $B\cap \nargmin f^\nu$ is nonempty for all $\nu$. With the refinement of tightness, we obtain the following consequences of epi-convergence (see \cite[Theorem 3.8]{RoysetWets.19a}), where
\[
\epsilon\mbox{-}\hspace{-0.05cm}\nargmin f = \{x\in \dom f ~|~ f(x) \leq \inf f + \epsilon\} ~~\mbox{for $f:\reals^n\to \Reals$ and $\epsilon \in \reals_+ = [0,\infty)$.}
\]

\begin{theorem}
\label{tEpiConvConsequences} {\rm (consequences of epi-convergence)}. For proper $f,f^\nu:\reals^n\to\Reals$, suppose that $f^\nu\eto f$. Then, the following
hold:
\begin{enumerate}[{\rm (a)}]

\item $\forall \big\{\epsilon^\nu\in \reals_+\to 0, \nu\in\nats\big\}$,
     $\nOutLim \big( \epsilon^\nu \mbox{-}\nargmin f^\nu\big) \subset \nargmin  f$.

\item $\nlimsup\,(\ninf  f^\nu) \leq \ninf  f$.

\item If $\big\{x^\nu\in \nargmin f^\nu, \nu\in N\big\}$ converges for some $N\in
    \cN_\infty^\grill$, then $\nlim_{\nu\in N} (\ninf f^{\nu}) = \ninf f$.

\item $\ninf  f^\nu \to \ninf f>-\infty \Longleftrightarrow \{f^\nu, \nu\in\nats\} \mbox{ is tight}$.

\item $\ninf  f^\nu \to \ninf f$ and $\epsilon>0$ $\Longrightarrow$
    $\nInnLim \big( \epsilon\mbox{-}\nargmin f^\nu\big) \supset
    \nargmin f$.

\item $\ninf  f^\nu \to \ninf f$  $\Longrightarrow$ $\exists \big\{\epsilon^\nu\in \reals_+\to 0, \nu\in\nats\big\}$ such that  $\nLim \big( \epsilon^\nu \mbox{-}\nargmin f^\nu\big) = \nargmin  f$.

\item $\exists \big\{\epsilon^\nu\in \reals_+\to 0, \nu\in\nats\big\}$ such that  $\nLim \big( \epsilon^\nu \mbox{-}\nargmin f^\nu\big) = \nargmin  f \neq \emptyset$  $\Longrightarrow$  $\ninf  f^\nu \to \ninf f>-\infty$.

\end{enumerate}
\end{theorem}
Part (a) formalizes the earlier assertion about minimizers. In light of the discussion prior to the theorem, (b,c) is as much as one can hope to say about minimum values without tightness. Part (f) supplements (a) by establishing set-convergence to $\nargmin f$ when the tolerance $\epsilon^\nu$ vanishes slowly enough.

\section{Variational Geometry and Subdifferentiability}\label{sec:vargeom}

Normal and tangent spaces are well-known from Differential Geometry as means to locally characterize sets defined by smooth manifolds. However, they fall short of the demand of variational problems where nonsmoothness is prevalent. Set-convergence supports the development of extensions of these concepts. This brings us to normal cones and tangent cones.

The {\it tangent cone}\footnote{A tangent cone is sometimes called a contingent cone.} to $C\subset\reals^n$ at one of its points $\bar x$ is defined as\footnote{For $C\subset \reals^n$, $\lambda \in \reals$, and $a\in \reals^n$: $\lambda C = \{\lambda x~|~x \in C\}$ (scalar multiplication) and $C+a = \{x+a~|~x\in C\}$ (addition).}
\[
T_C(\bar x) = \nOutLim \nu(C-\bar x), \mbox{ with } w\in T_C(\bar x) \mbox{ being a {\it tangent vector} to } C \mbox{ at } \bar x.
\]
Geometrically, $C-\bar x$ is a translation of $C$ to the origin such that $0$ occupies the same relative position in $C-\bar x$ as $\bar x$ in $C$ and $T_C(\bar x)$ is the outer limit of this translation, ``magnified'' by $\nu$; see Figure \ref{fig:deftangent}.

\drawing{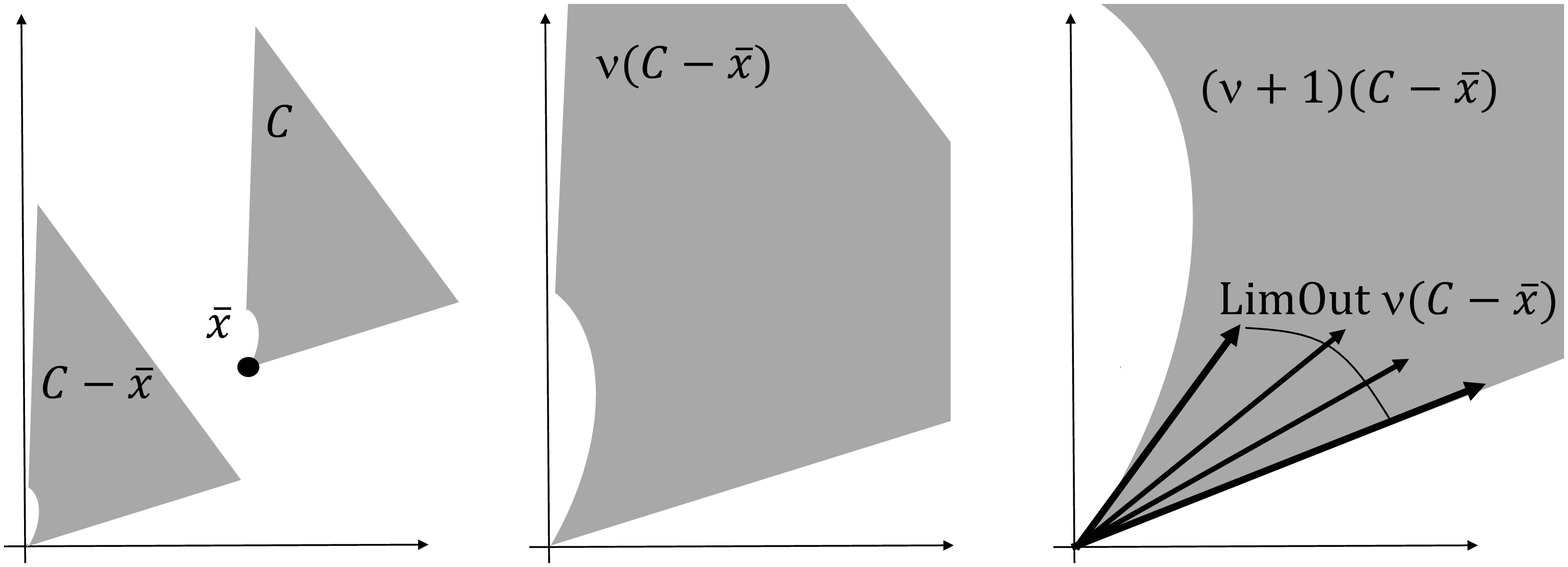}{4.2in}
   {Translation (left), magnification (middle), further magnification and tangent cone (right).} {fig:deftangent}

Figure \ref{fig:tangentcone} illustrates typical situations. In the left
portion, $C$ has a ``smooth'' boundary at $\bar x$ and $T_C(\bar x)$ is simply the halfspace formed by the usual tangent space of the boundary of $C$ at
$\bar x$. When $\bar x$ is at a ``corner,'' as seen in the middle portion of
the figure, the tangent cone narrows. At the ``inward kink'' of the heart-shaped set to the right, the tangent cone fans out in all directions except some pointing upward. In this figure and elsewhere, tangent vectors are often visualized
after being translated by $\bar x$, i.e., the beginning of an arrow visualizing a tangent vector $w$ is moved from the
origin, where it always will be, to $\bar x$ so the illustrations are actually of $\bar x + w$ and not $w$.

\drawing{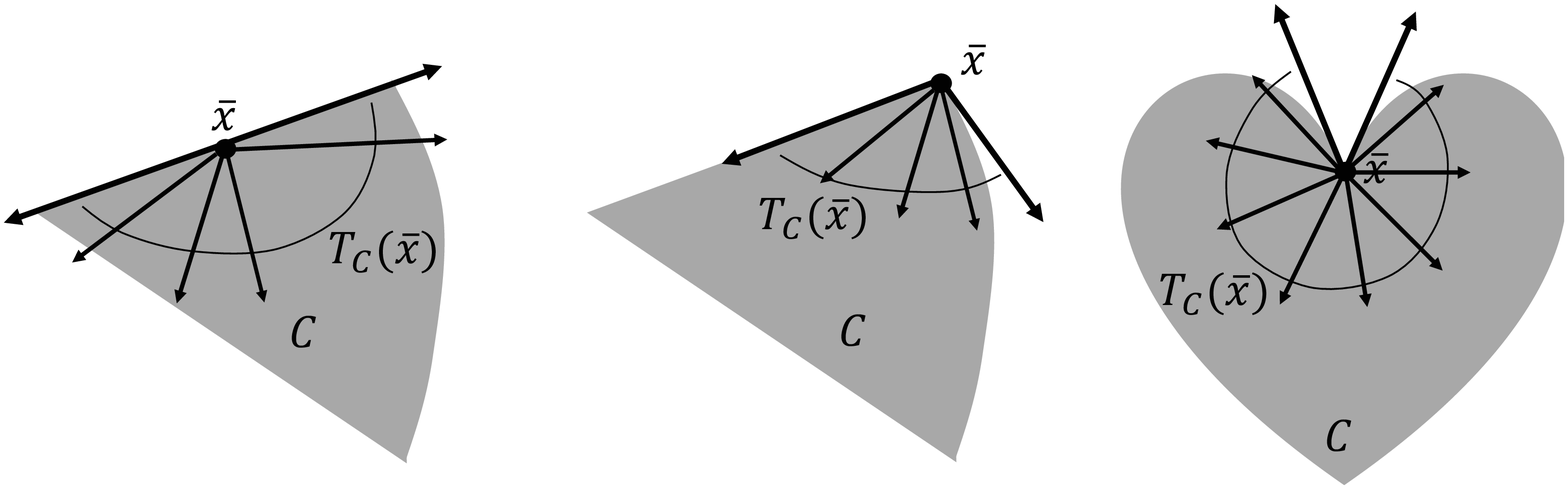}{4.5in}
   {Tangent vectors and tangent cones translated to $\bar x$.} {fig:tangentcone}

For $\bar x \in C\subset \reals^n$, $v\in\reals^n$ is a {\it regular normal vector} to $C$ at $\bar x$ if
\[
\langle v, w\rangle \leq 0 \mbox{ for all } w\in T_C(\bar x).
\]
The set of all such regular normal vectors is $\widehat N_C(\bar x)$. The {\it normal cone} to $C$ at $\bar x$ is defined as
\[
N_C(\bar x) = \bigcup_{x^\nu\in C\to \bar x} \nOutLim \widehat N_C(x^\nu), ~\mbox{ with } v\in N_C(\bar x) \mbox{ being a {\it normal vector} to } C \mbox{ at } \bar x.
\]
We recall that $\langle v,w\rangle = \|v\|_2 \|w\|_2 \cos \alpha$, where
$\alpha$ is the angle between the vectors $v$ and $w$. Thus, a regular normal vector forms an angle of between 90 and 270 degrees with every tangent vector and thereby extends the ``perpendicular relationship'' of normal and tangent spaces from Differential Geometry. Figure \ref{fig:normalcone2} illustrates the normal cones corresponding to the tangent cones of Figure \ref{fig:tangentcone}.

\drawing{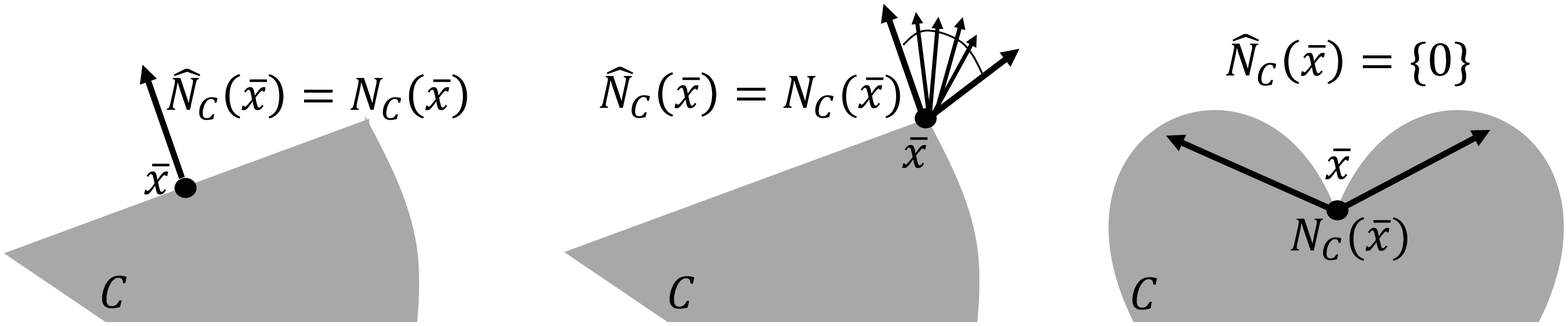}{4.5in}
   {Normal vectors and normal cones translated to $\bar x$.} {fig:normalcone2}

On their own, regular normal vectors fail to provide a solid basis for the treatment of ``complicated'' sets as seen in the right portion of Figure \ref{fig:normalcone2}. At the inward kink in this heart-shaped set, $\widehat N_C(\bar x)=\{0\}$ since $T_C(\bar x)$ fans out in nearly all directions and it becomes impossible to find a nonzero vector that forms an angle of at least 90 degrees with all tangent vectors. We obtain a robust notion of ``normality'' by considering the union of the outer limits of $\widehat N_C(x^\nu)$ for all $x^\nu\in C\to \bar x$. This leads to vectors in $N_C(\bar x)$ that aren't regular normals. Figure \ref{fig:normalcone5} shows $v^\nu \in \widehat N_C(x^\nu)$ at points $x^\nu\in C$ approaching $\bar x$ from the left. In the limit, these regular normal vectors give rise to the normal vectors pointing northeast in the right-most portion of the figure. Likewise, $x^\nu\in C$ approaching $\bar x$ from the right result in the normal vectors at $\bar x$ pointing northwest. It might be counterintuitive that the normal vectors at $\bar x$ points {\it into} the set, but this is indeed representative of the local geometry at this point and allows us to treat such points in a meaningful way. In general, a normal cone is a closed set.

\drawing{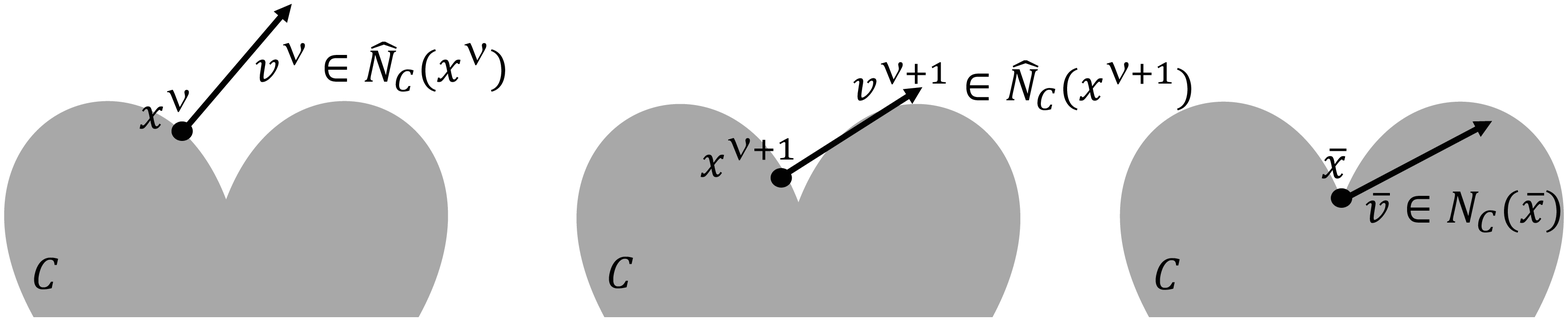}{4.5in}
   {Normal vectors as limits of regular normal vectors.} {fig:normalcone5}

\begin{example}{\rm (optimality conditions)}.\label{eOptCond} For a smooth\footnote{A function is smooth if its gradient is continuous.} function $f:\reals^n\to \reals$ and a closed set $C\subset\reals^n$, consider the problem of minimizing $f$ over $C$. Then,
\[
x^\star \mbox{ is a local minimizer } ~~\Longrightarrow~~ -\nabla f(x^\star) \in N_C(x^\star),
\]
which can be seen, for example, by invoking \cite[Exercise 10.10]{VaAn}. Suppose that a numerical procedure generates the sequence $\{x^\nu \in C, \nu\in\nats\}$ with $\dist_2(-\nabla f(x^\nu), N_C(x^\nu)) \leq \epsilon^\nu$ for some tolerance $\epsilon^\nu \downto 0$. Then, the continuity property built into normal cones ensures that any cluster point $\bar x$ of the sequence satisfies $-\nabla f(\bar x) \in N_C(\bar x)$.
\end{example}
\state Detail. Suppose that $x^\nu\Nto \bar x$ for some $N\in \cN_\infty^\grill$.  By continuity, $\nabla f(x^\nu) \Nto \nabla f(\bar x)$. Let $v^\nu \in \nargmin\{ \|v - \nabla f(x^\nu)\|_2~ | -v\in N_C(x^\nu) \}$. Then, $\|v^\nu - \nabla f(x^\nu)\|_2\leq \epsilon^\nu$, which implies that $v^\nu \Nto \nabla f(\bar x)$. Since $-v^\nu\in N_C(x^\nu) \Nto -\nabla f(\bar x)$, we must have $-\nabla f(\bar x) \in \nOutLim_{\nu\in N} N_C(x^\nu)$.  From the definition of normal cones, $\nOutLim_{\nu\in N} N_C(x^\nu) \subset N_C(\bar x)$ because $x^\nu \Nto \bar x$ and the assertion follows.\eop

We next turn the attention to subgradients, which are defined in terms of normal cones of epigraphs.
For $f:\reals^n\to \Reals$ and $\bar x$ with $f(\bar x)$ finite, $v$ is a {\it subgradient}\footnote{These subgradients are also called Mordukhovich or limiting subgradients, to distinguish them from other kinds, but we omit such qualifications because, at least in finite dimensions, they are the central ones; see \cite[Chapter 8]{VaAn}. For convex functions, these subgradients coincide with those from convex analysis.} of $f$ at $\bar x$ if
\[
(v,-1)\in N_{\epi f} \big(\bar x, f(\bar x)\big).
\]
The set of all subgradients of $f$ at $\bar x$ is denoted by $\partial f(\bar x)$.

As an introductory example, let $f(x) =
x^2$ for $x\in\reals$. Then, $N_{\epi f}(1, f(1))$ $=$ $\{(2,-1)y~|~y\geq 0\}$. Thus, $(v,-1) \in
N_{\epi f}(1, f(1))$ if and only if $v = 2$ and $\partial f(1) = \{2\}$ as expected. Figure \ref{fig:epideriv4} illustrates the normal cones of an epigraph at two points. Roughly, it seems like $N_{\epi f} (\bar x, f(\bar x))$ contains vectors of the form $\lambda(v,-1)$ for $\lambda \geq 0$ and $v \in [1/2, 2]$; these are the subgradients of $f$ at $\bar x$. Note that $2$ is the slope of $f$ to the right of $\bar x$ and $1/2$ is the slope to the left. At $(x', f(x'))$, normal vectors appear to be of the form $\lambda(-4, -1)$ and $\lambda(-1/2,-1)$ for $\lambda\geq 0$. Then, $v = -4$ and $v = -1/2$ are subgradients and in fact the only ones because the normal cone at $(x', f(x'))$ consists of only two rays.

\drawing{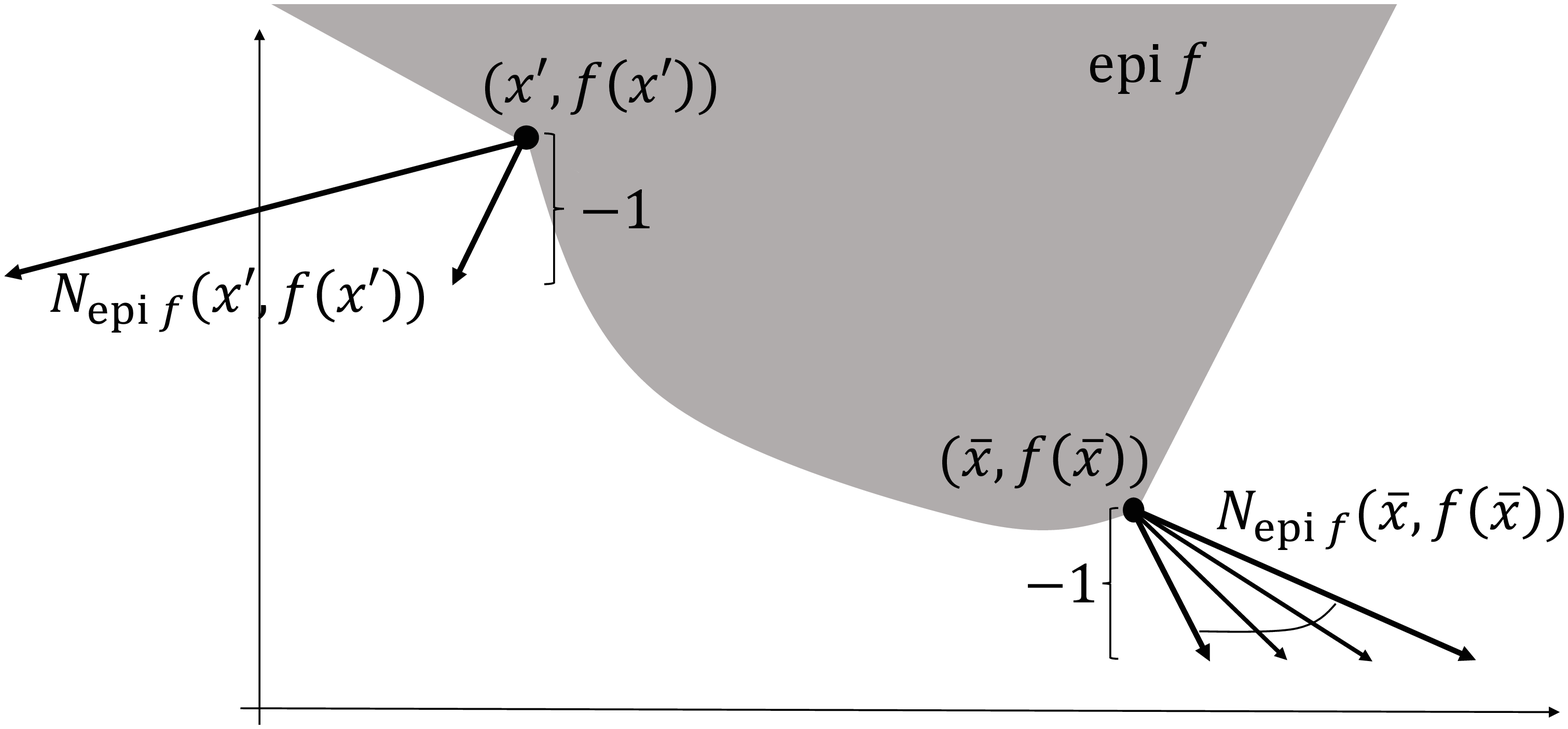}{3.5in} {Normal cones of $\epi f$.} {fig:epideriv4}

\begin{example}{\rm (the Fermat rule).}\label{eFermat}
For $f:\reals^n\to \Reals$ and $x^\star$ with $f(x^\star)>-\infty$,
\[
x^\star \mbox{ is a local minimizer of } f ~~\Longrightarrow~~ 0 \in \partial f(x^\star),
\]
which is the celebrated Fermat rule; see \cite[Theorem 8.15]{VaAn}. If an algorithm generates a sequence $\{x^\nu, \nu\in\nats\}$ with $\dist_2(0, \partial f(x^\nu)) \leq \epsilon^\nu$ for some tolerance $\epsilon^\nu \downto 0$, $\bar x$ is a cluster point, i.e., $x^\nu\Nto \bar x$ for some $N\in \cN_\infty^\grill$, and $f(x^\nu) \Nto f(\bar x)\in\reals$, then $\bar x$ satisfies the Fermat Rule, i.e, $0 \in \partial f(\bar x)$.
\end{example}
\state Detail. In light of Example \ref{eOptCond} and the definition of subgradients,
$\nOutLim_{\nu\in N} \partial f(x^\nu) \subset \partial f(\bar x)$; see also \cite[Proposition 8.7]{VaAn}. Arguments similar to the ones in Example \ref{eOptCond} establish the assertion.\eop

These concepts are supplemented by subderivatives of functions (akin to directional derivatives) which are defined in terms of outer limits of epigraphs of certain difference quotients, with further extensions leading to notions of differentiation of set-valued mappings again utilizing set-convergence as a primary tool; see \cite[Chapter 8]{VaAn}. We refer to the commentary sections of \cite[Chapters 6 and 8]{VaAn} for the historical development of these concepts.

\section{Generalized Equations and Graphical Convergence}\label{sec:geneq}

While we express a system of equations $F(x) = 0$ in terms of a mapping $F:\reals^n\to \reals^m$, a generalized equation is given by a {\it set-valued mapping} $S:\reals^n\tto \reals^m$. The distinguishing feature of $S$ relative to $F$ is that we now allow for an input $x$ to produce an output $S(x)$ that is a subset of $\reals^m$ and not necessarily a single point. In fact, we even allow the subset to be empty so that the domain
\[
\dom S = \big\{x\in\reals^n~\big|~S(x) \neq\emptyset \big\}
\]
could be a strict subset of $\reals^n$. Of course, $F$ defines a particular set-valued mapping with $S(x) = \{F(x)\}$ for all $x\in \reals^n$ with $\dom S = \reals^n$.  The {\it graph}\index{graph} of a set-valued mapping $S:\reals^n\tto\reals^m$ is the set
\[
\gph S = \big\{(x,y)\in\reals^n\times \reals^m~\big|~y\in S(x)\big\},
\]
which describes $S$ fully; see Figure \ref{fig_setval3}. Given $\bar y\in\reals^m$ and $S:\reals^n\tto\reals^m$, a {\it generalized equation} is then
\[
\bar y \in S(x), ~~\mbox{ with set of solutions } S^{-1}(\bar y) = \big\{ x\in \reals^n  ~\big|~ \bar y \in S(x)\big\}.
\]

\drawing{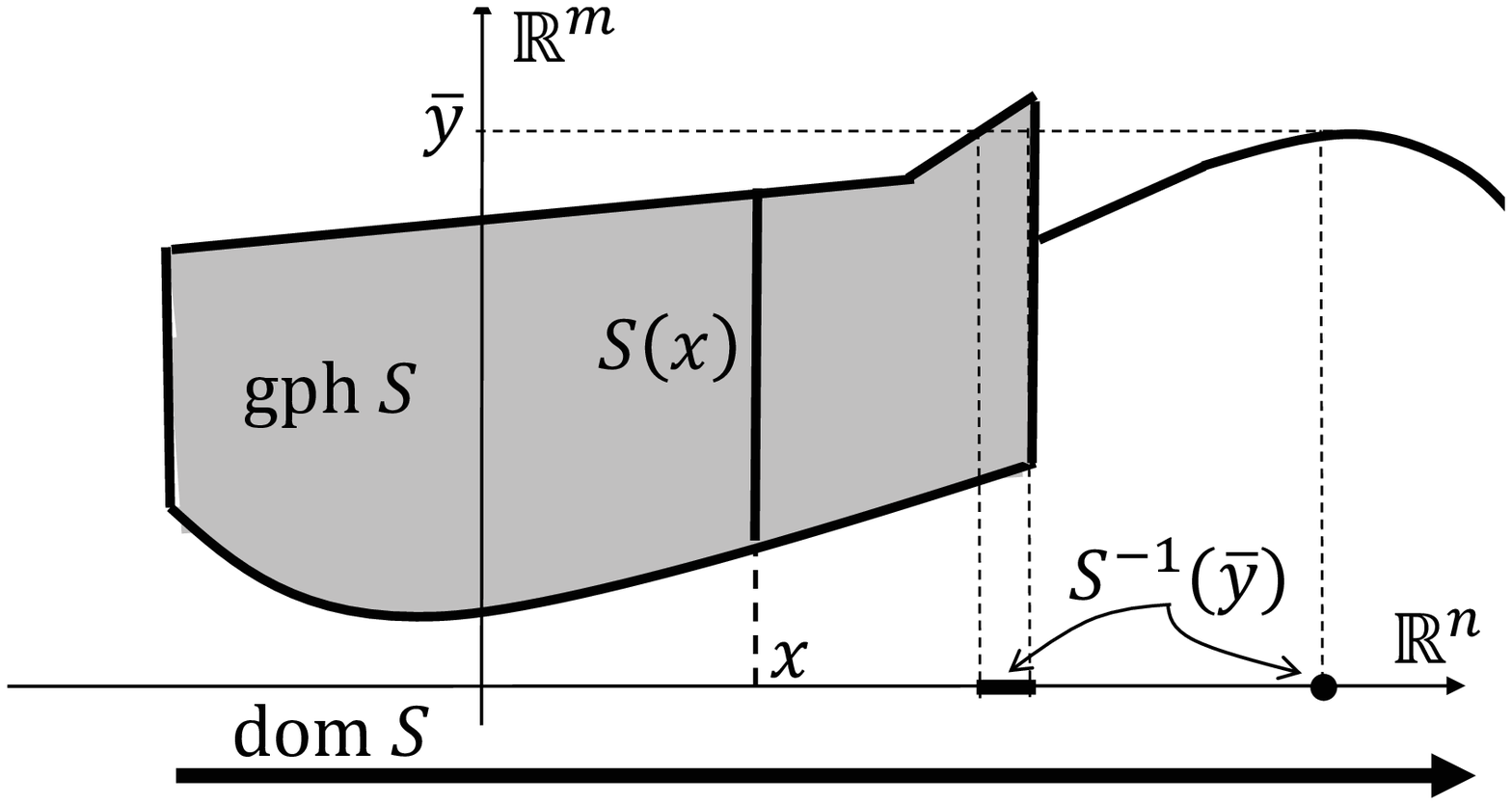}{2.7in} {The solution set $S^{-1}(\bar y)$ to the generalized equation $\bar y \in S(x)$.} {fig_setval3}

The sets of subgradients of $f:\reals^n\to \Reals$ can be viewed as a set-valued mapping $S:\reals^n\tto \reals^n$ with $S(x) = \partial f(x)$ if $f(x)$ is finite and $S(x) = \emptyset$ otherwise. Likewise, the normal cones of $C\subset\reals^n$ define a set-valued mapping $S:\reals^n\tto \reals^n$ with $S(x) = N_C(x)$ if $x\in C$ and $S(x) = \emptyset$ otherwise. The facts about $\nOutLim N_C(x^\nu)\subset N_C(\bar x)$ and $\nOutLim \partial f(x^\nu)\subset \partial f(\bar x)$ that we utilize in Examples \ref{eOptCond} and \ref{eFermat} can be viewed as continuity properties of the corresponding set-valued mappings.

Generally, a set-valued mapping $S:\reals^n\tto\reals^m$ is {\it outer semicontinuous} (osc)\index{outer semicontinuous} at $\bar x\in\reals^n$ if
\[
\bigcup_{x^\nu\to \bar x} \nOutLim S(x^\nu) \subset S(\bar x).
\]
It's osc if this holds for all $\bar x\in\reals^n$. Figure \ref{fig_setval3} furnishes an example.

\begin{proposition}{\rm (characterization of osc).}\label{pCharOSC}
For $S:\reals^n\tto \reals^m$,
\[
\gph S \mbox{ is closed } ~~\Longleftrightarrow~~ S \mbox{ is osc } ~~\Longleftrightarrow~~ S^{-1} \mbox{ is osc}.
\]
\end{proposition}
\state Proof. A closed $\gph S$ means that every sequence $(x^\nu,v^\nu)\in \gph S\to (\bar x, \bar v)$ has $(\bar x, \bar v) \in \gph S$. But, this is equivalent to $x^\nu\to \bar x$, $v^\nu\to \bar v$, and $v^\nu\in S(x^\nu)$ implying $\bar v\in S(\bar x)$, which is guaranteed by osc. The claim about the inverse mapping follows similarly because $v \in S(x)$ if and only if $x\in S^{-1}(v)$.\eop

\begin{example}{\rm (osc of feasible-set mapping).}\label{eFeasiblemap} For a closed $X\subset \reals^n$ and continuous $f_i,g_i:\reals^d\times \reals^n\to \reals$, consider the set-valued mapping $S:\reals^d\tto \reals^n$ given by
\[
S(u) = \big\{ x\in X  ~\big|~ f_i(u, x) = 0, i=1, \dots, m; ~~ g_i(u, x) \leq 0, i=1, \dots, q\big\},
\]
which arises in stability studies of feasible sets under perturbations. In fact, $S$ is osc.
\end{example}
\state Detail. In view of Proposition \ref{pCharOSC}, it suffices to establish the closedness of the graph
\[
\gph S = \big\{(u,x)\in\reals^d\times\reals^n~\big|~ x\in S(u)\big\} =  \bigcap_{i=1}^m  \big\{ (u,x) \big| f_i(u, x) = 0\big\}\bigcap_{i=1}^q  \big\{ (u,x) \big| g_i(u, x) \leq 0\big\} \bigcap \big( \reals^d \times X \big).
\]
Since $X$ is closed and the functions involved are continuous, this is an intersection of closed sets and is therefore closed.\eop

A set-valued mapping $S:\reals^n\tto\reals^m$ is {\it inner semicontinuous} (isc) at $\bar x\in\reals^n$ if
\[
\bigcap_{x^\nu\to \bar x} \nInnLim S(x^\nu) \supset S(\bar x).
\]
It's isc if this holds for all $\bar x\in\reals^n$. Moreover, $S$ is {\it continuous} (at $\bar x$) if it's both osc and isc (at $\bar x$), in which case $S(x^\nu)\sto S(\bar x)$ for all $x^\nu\to \bar x$.

While osc is common in applications, isc occurs less frequently. As an example, if $S:\reals^n\tto \reals^m$ has a convex graph, then $S$ is isc at any point $\bar x \in \nt (\dom S)$; see \cite[Theorem 5.9b]{VaAn}.

Applications often bring us generalized equations that can't be solved by existing algorithms, at least not in reasonable time,  and it becomes necessary to consider approximations. Approximations may also be introduced artificially in sensitivity analysis. Regardless of the circumstances, we would like to examine whether replacing $S:\reals^n\tto\reals^m$ by an alternative $S^\nu:\reals^n\tto \reals^m$ results in significantly different solutions of the corresponding generalized equations. It's immediately clear that if an algorithm computes $\{x^\nu, \nu\in\nats\}$ with the property that $\dist_2(\bar y, S^\nu(x^\nu))\leq \epsilon^\nu$ for some vanishing tolerances $\epsilon^\nu\geq 0$, then every cluster point $\bar x$ of the sequence satisfies $\bar y\in S(\bar x)$ provided that
\[
\nOutLim (\gph S^\nu) \subset \gph S.
\]
To see this, note that the vanishing $\dist_2(\bar y, S^\nu(x^\nu))$ guarantees that there is $y^\nu\in S^\nu(x^\nu)$ such that $y^\nu\to \bar y$. Since $(x^\nu,y^\nu)\in \gph S^\nu$, it follows by the definition of outer limits that $(\bar x, \bar y) \in \nOutLim (\gph S^\nu)$ and then also $(\bar x, \bar y) \in \gph S$ under the stated assumption. Consequently, the algorithm is justified.

\begin{example}{\rm (smoothing of complementarity problems).}\label{eCompl}
For a smooth\footnote{A mapping is smooth if its Jacobian is continuous.} mapping $F:\reals^n\to \reals^n$, the complementarity problem
\[
x\geq 0, ~~~F(x)\geq 0, ~~~ \big\langle F(x), x\big\rangle = 0
\]
corresponds to the generalized equation $0 \in F(x) + N_C(x)$, with $C = \reals_+^n$. This, in turn, is equivalent to the normal map condition 
\[
0 = F\big( \prj_C(z) \big) + z - \prj_C(z),
\]
where the projection $\prj_C(z) = (\max\{0,z_1\}, \dots, \max\{0,z_n\})$; see \cite{Robinson.92}. Although just ``regular'' equations, the normal map condition involves nonsmooth functions and the usual Newton's Method doesn't apply. However, a smooth approximation can be brought in.
\end{example}
\state Detail.  For $\theta^\nu>0$, $\phi^\nu(\alpha) = \frac{1}{\theta^\nu}\ln\big(1 + e^{ \alpha \theta^\nu}\big)$ furnishes a smooth approximation of $\alpha\mapsto \max\{0,\alpha\}$ with error
\[
0\leq  \phi^\nu(\alpha) - \max\{0,\alpha\}   \leq \frac{\ln 2}{\theta^\nu} ~\mbox{ for all } \alpha\in\reals.
\]
Consequently, the projection mapping can be approximated by the smooth mapping $\Phi^\nu:\reals^n\to \reals$ given by $\Phi^\nu(z) = (\phi^\nu(z_1), \dots, \phi^\nu(z_n))$. This results in the approximating equations
\[
0 = F\big( \Phi^\nu(z) \big) + z - \Phi^\nu(z),
\]
which can be solved using Newton's Method. We justify this approach by viewing the equations in terms of the set-valued mappings $S,S^\nu:\reals^n\tto\reals^n$ given by $S(z) = \{F( \prj_C(z) ) + z - \prj_C(z)\}$ and $S^\nu(z) = \{F( \Phi^\nu(z) ) + z - \Phi^\nu(z)\}$ for all $z\in\reals^n$. Clearly, the normal map condition is then $0\in S(z)$. In view of the discussion prior to the example, we only need to establish that $\nOutLim (\gph S^\nu) \subset \gph S$.

Suppose that $(\bar z, \bar v)\in \nOutLim (\gph S^\nu)$. Then, there are $N\in \cN_\infty^\grill$, $z^\nu\Nto \bar z$, and $v^\nu\Nto \bar v$ such that $v^\nu\in S^\nu(z^\nu)$ for $\nu\in N$. Moreover,
\begin{align*}
\big\|F\big( \prj_C(\bar z) \big) + \bar z - \prj_C(\bar z) - \bar v\big\|_2 & \leq  \big\|F\big( \prj_C(\bar z) \big) + \bar z - \prj_C(\bar z) - F\big( \Phi^\nu(z^\nu) \big) - z^\nu + \Phi^\nu(z^\nu)\big\|_2\\
&~~~~+ \big\| F\big( \Phi^\nu(z^\nu) \big) + z^\nu - \Phi^\nu(z^\nu) - \bar v\big\|_2.
\end{align*}
Since $F( \Phi^\nu(z^\nu) ) + z^\nu - \Phi^\nu(z^\nu)= v^\nu\Nto \bar v$, the second term on the right-hand side vanishes. Also,
\begin{align*}
&\big\|F\big( \prj_C(\bar z) \big) + \bar z - \prj_C(\bar z) - F\big( \Phi^\nu(z^\nu) \big) - z^\nu + \Phi^\nu(z^\nu)\big\|_2\\
\leq &\big\|F\big( \prj_C(\bar z) \big) - F\big( \Phi^\nu(z^\nu) \big) \big\|_2 + \|\bar z - z^\nu\|_2 + \big\|\prj_C(\bar z) - \Phi^\nu(z^\nu)\big\|_2,
\end{align*}
which vanishes as long as $\theta^\nu\to \infty$ because then $\phi^\nu(z_j^\nu)\Nto \max\{0,\bar z_j\}$ for all $j=1, \dots, n$. Consequently, $F( \prj_C(\bar z) ) + \bar z - \prj_C(\bar z) = \bar v$ and $(\bar z, \bar v) \in \gph S$.\eop

With the pioneering work of Spagnolo \cite{Spagnolo.76}, Attouch \cite{Attouch.77}, De Giorgi \cite{Degiorgi.77}, and Moreau \cite{Moreau.78}, set-convergence of graphs of set-valued mappings, which became known as graphical convergence, and its importance for generalized equations came to the foreground. For $S^\nu, S:\reals^n\tto \reals^m$, we say that $S^\nu$ {\it converges graphically} to $S$, written $S^\nu\gto S$, when
\[
\gph S^\nu \sto \gph S.
\]
In light of the definitions, this takes place when the condition $\nOutLim (\gph S^\nu) \subset \gph S$ is supplemented by $\nInnLim (\gph S^\nu) \supset \gph S$; see Figure \ref{fig:graphconv}.

Since $v\in S(x)$  if and only if $x\in S^{-1}(v)$, graphical convergence of set-valued mappings corresponds to graphical convergence of the corresponding inverse mappings:
\[
S^\nu \gto S ~~\Longleftrightarrow  ~~ (S^\nu)^{-1} \gto S^{-1}.
\]

\drawing{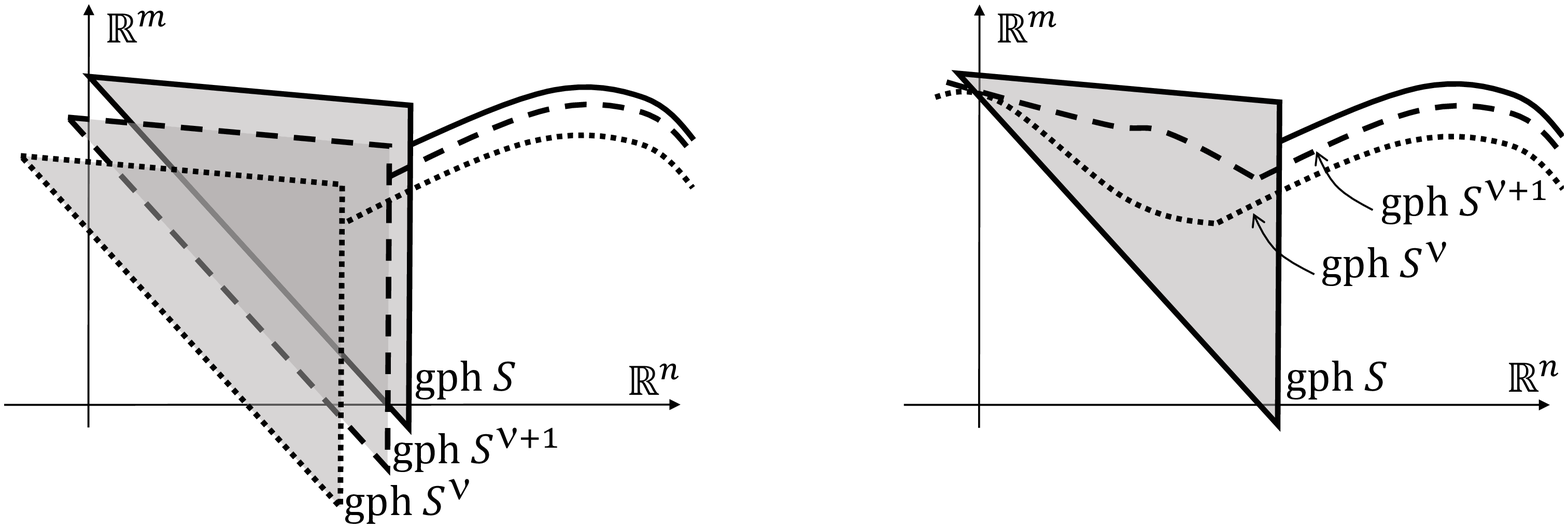}{5.0in} {On the left, $\gph S^\nu \sto \gph S$, but on the right only $\nOutLim (\gph S^\nu) \subset \gph S$ holds.} {fig:graphconv}

Under graphical convergence, conclusions about solutions of approximating generalized equations can be strengthen beyond the discussion prior to Example \ref{eCompl}. For $S:\reals^n\tto\reals^m$, $\bar y\in \reals$, and $\epsilon\geq 0$, the {\it set of $\epsilon$-solutions} to the generalized equation $\bar y \in S(x)$ is defined as
\[
S^{-1}\big( \ball(\bar y, \epsilon)  \big) = \bigcup_{y\in \ball(\bar y, \epsilon)} S^{-1}(y),
\]
where the ball $\ball(\bar y, \epsilon) = \{y\in \reals^m~|~\|y-\bar y\|\leq \epsilon\}$ is given in terms of any norm $\|\cdot\|$ on $\reals^m$.

\begin{theorem}\label{tgraphConv}{\rm (consequences of graphical convergence).} For set-valued mappings $S,S^\nu:\reals^n\tto\reals^m$, suppose that $S^\nu \gto S$ and $\bar y\in \reals^m$. Then, the following hold:
\begin{enumerate}[{\rm (a)}]
\item $\forall \{\epsilon^\nu\in \reals_+\to 0, \nu\in\nats\}$, $\nOutLim (S^\nu)^{-1}\big(\ball(\bar y, \epsilon^\nu)\big) \subset S^{-1}(\bar y)$.

\item $\exists \{\epsilon^\nu\in \reals_+\to 0, \nu\in\nats\}$,  with $\nInnLim (S^\nu)^{-1}\big(
   \ball(\bar y, \epsilon^\nu) \big) \supset S^{-1}(\bar y)$.
\end{enumerate}
\end{theorem}
\state Proof. For (a), we only utilize $\nOutLim (\gph (S^\nu)^{-1}) \subset \gph S^{-1}$: Let $\bar x \in \nOutLim (S^\nu)^{-1}(\ball(\bar y, \epsilon^\nu))$. Then, there exist $N\in \cN_\infty^\grill$ and $x^\nu\Nto \bar x$ such that $x^\nu \in (S^\nu)^{-1}(\ball(\bar y, \epsilon^\nu))$ and, for some $y^\nu\in \ball(\bar y, \epsilon^\nu)$, $(y^\nu,x^\nu) \in \gph (S^\nu)^{-1}$. Since cluster points of such sequences are contained in $\gph S^{-1}$, $\bar x\in S^{-1}(\bar y)$.

For (b), we use $\nInnLim (\gph (S^\nu)^{-1}) \supset \gph S^{-1}$: Let $\bar x\in S^{-1}(\bar y)$. Then, $(\bar y, \bar x) \in \nInnLim (\gph (S^\nu)^{-1})$ and there exist $y^\nu\to \bar y$ and $x^\nu\to \bar x$ such that $x^\nu\in (S^\nu)^{-1}(y^\nu)$. We seek to establish that $\bar x \in \nInnLim (S^\nu)^{-1}(\ball(\bar y, \epsilon^\nu))$ for some vanishing $\epsilon^\nu$ and need to construct $\bar x^\nu\to \bar x$, $\bar y^\nu\to \bar y$, and $\epsilon^\nu\to 0$ such that $\bar x^\nu\in (S^\nu)^{-1}(\bar y^\nu)$ and $\bar y^\nu \in \ball(\bar y, \epsilon^\nu)$. We see that $\bar x^\nu = x^\nu$, $\bar y^\nu = y^\nu$, and $\epsilon^\nu = \|y^\nu-\bar y\|$ suffices.\eop

The theorem establishes that cluster points of near-solutions of approximating generalized equations are indeed solutions of the actual generalized equation as long as the tolerance vanishes. Moreover, every such solution can be approached by near-solutions of the approximating generalized equations provided that the tolerance vanishes sufficiently slowly.

\begin{example}{\rm (homotopy methods for generalized equations).}\label{eHomo}
For an osc set-valued mapping $S:\reals^n\tto \reals^n$ and $\bar y\in \reals^n$, a homotopy method for solving the generalized equation $\bar y \in S(x)$ involves solving the approximating generalized equations
\[
\bar y^\nu \in S^\nu(x) = (1-\lambda^\nu) S(x) + \lambda^\nu x
\]
for $\lambda^\nu \in (0,1]\to 0$ and $\bar y^\nu\to \bar y$. Trivially, $\bar y^\nu$ solves this equation when $\lambda^\nu = 1$ and this provides a starting point for iterations involving successively smaller $\lambda^\nu$. The hope is that a gradual change in $\lambda^\nu$ results in approximating generalized equations that can be solve much quicker than the actual problem. In lieu of Theorem \ref{tgraphConv}, we justify the method by establishing that $S^\nu\gto S$.
\end{example}
\state Detail. Let $(\hat x,\hat y) \in \nOutLim (\gph S^\nu)$. Then, there are $N\in \cN_\infty^\grill$, $x^\nu \Nto \hat x$, and $y^\nu \Nto \hat y$ with $y^\nu\in S^\nu(x^\nu)$. Consequently, $y^\nu = (1-\lambda^\nu)v^\nu + \lambda^\nu x^\nu$ for some $v^\nu\in S(x^\nu)$. Since $v^\nu = (1-\lambda^\nu)^{-1}(y^\nu - \lambda^\nu x^\nu) \Nto \hat y$ and $S$ is osc, $\hat y \in S(\hat x)$. This establishes that $\nOutLim (\gph S^\nu) \subset \gph S$.

Next, take $(\hat x, \hat y) \in \gph S$ and construct $x^\nu = \hat x$ and $y^\nu = (1-\lambda^\nu )\hat y + \lambda^\nu \hat x$. Then, $(x^\nu,y^\nu) \in \gph S^\nu$, $x^\nu \to \hat x$, and $y^\nu \to \hat y$, which means that $(\hat x, \hat y) \in \nInnLim (\gph S^\nu)$. This establishes that $\nInnLim (\gph S^\nu) \supset \gph S$ and then also $S^\nu\gto S$.

Figure \ref{fig_homotopy} illustrates the approach when applied to $0 \in S(x) = \{x + \sin x + 1\}$, which is associated with the challenge that the derivative of $x\mapsto x + \sin x + 1$ vanishes for $x = \pi + 2\pi k$, $k \in \{\dots, -2, -1, 0, 1, 2, \dots\}$ and Newton's Method breaks down at these points. However, the approximation $S^\nu(x) = \{(1-\lambda^\nu)(x + \sin x + 1) + \lambda^\nu x\}$ doesn't have that problem for $\lambda^\nu>0$. The figure shows the graph of $S$ (solid line) and those of $S^\nu$ for $\lambda^\nu = 1, 0.5, 0.1$.\eop

\drawing{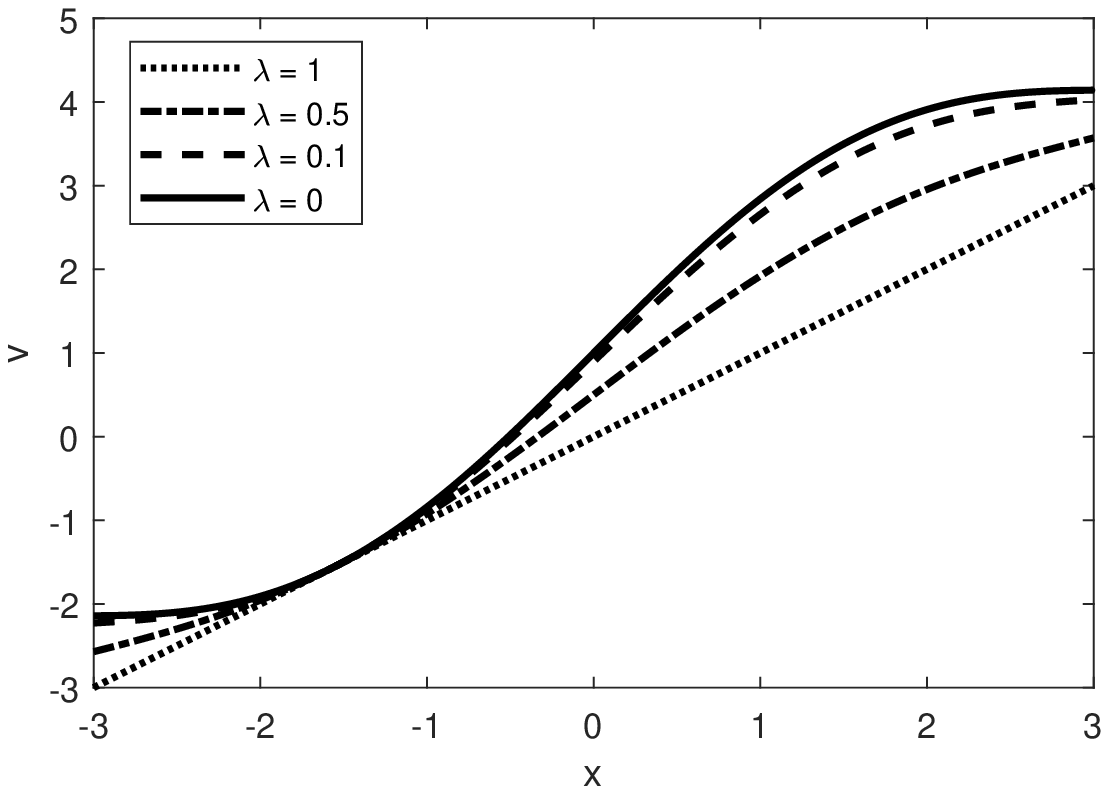}{3.6in} {Graph of $S$ (solid line) in Example \ref{eHomo} and those of its graphically converging approximations.} {fig_homotopy}

In the convex case, we have a useful connection  due to Attouch \cite{Attouch.77} between epi-convergence of functions and graphical convergence of their sets of subgradients viewed as set-valued mappings.

\begin{proposition}{\rm (convergence for convex functions).}
For proper, lsc, convex functions $f^\nu,f:\reals^n\to \Reals$,
\[
  f^\nu \eto f ~\Longleftrightarrow~ \partial f^\nu \gto \partial f \mbox{ and } f^\nu(x^\nu) \to f(\bar x)\mbox{ for some } v^\nu\to \bar v, x^\nu \to \bar x \mbox{ with } v^\nu \in \partial f^\nu(x^\nu) \mbox{ and } \bar v \in \partial f(\bar x).
\]
\end{proposition}

\section{Quantification of Set-Convergence}\label{sec:quant}

The rate of convergence of approximating solutions to an actual one relates to the rate with which the corresponding set-convergence takes place. This, in turn, requires us to clarify
the meaning of distances between two sets. To allow for more flexibility, we define parallel to \eqref{eqn:distdef}
the point-to-set distance $\dist(\bar x,C) = \ninf_{x\in C} \|x-\bar x\|$ for a nonempty $C\subset\reals^n$ using any norm $\|\cdot\|$ on $\reals^n$ and $\dist(\bar x, \emptyset) = \infty$.

The {\it excess} of $C$ over $D$ is defined as
\[
\exs(C; D) =\begin{cases}
\nsup_{x\in C} \dist(x,D) & \mbox{ if } C \neq\emptyset, D\neq\emptyset\\
\infty & \mbox{ if } C \neq\emptyset, D=\emptyset\\
0 & \mbox{ otherwise}.
\end{cases}
\]
Figure \ref{fig_excess} shows that $\exs(C;D)$ can be different than $\exs(D;C)$ and it becomes natural to look at the sum of the two, as proposed by Pompeiu in 1905, or the larger of the two as promoted by Hausdorff in 1927. The latter has become known as the Hausdorff distance. In our setting with unbounded sets commonly occurring, either of these could be infinity.
Thus, we adopt a truncation that can be traced back to Walkup and Wets \cite{WalkupWets.67} and Mosco \cite{Mosco.69} in their work on convex cones and convex sets, respectively. Its systematic study started with the work of Attouch and Wets \cite{AttouchWets.91}.

\drawing{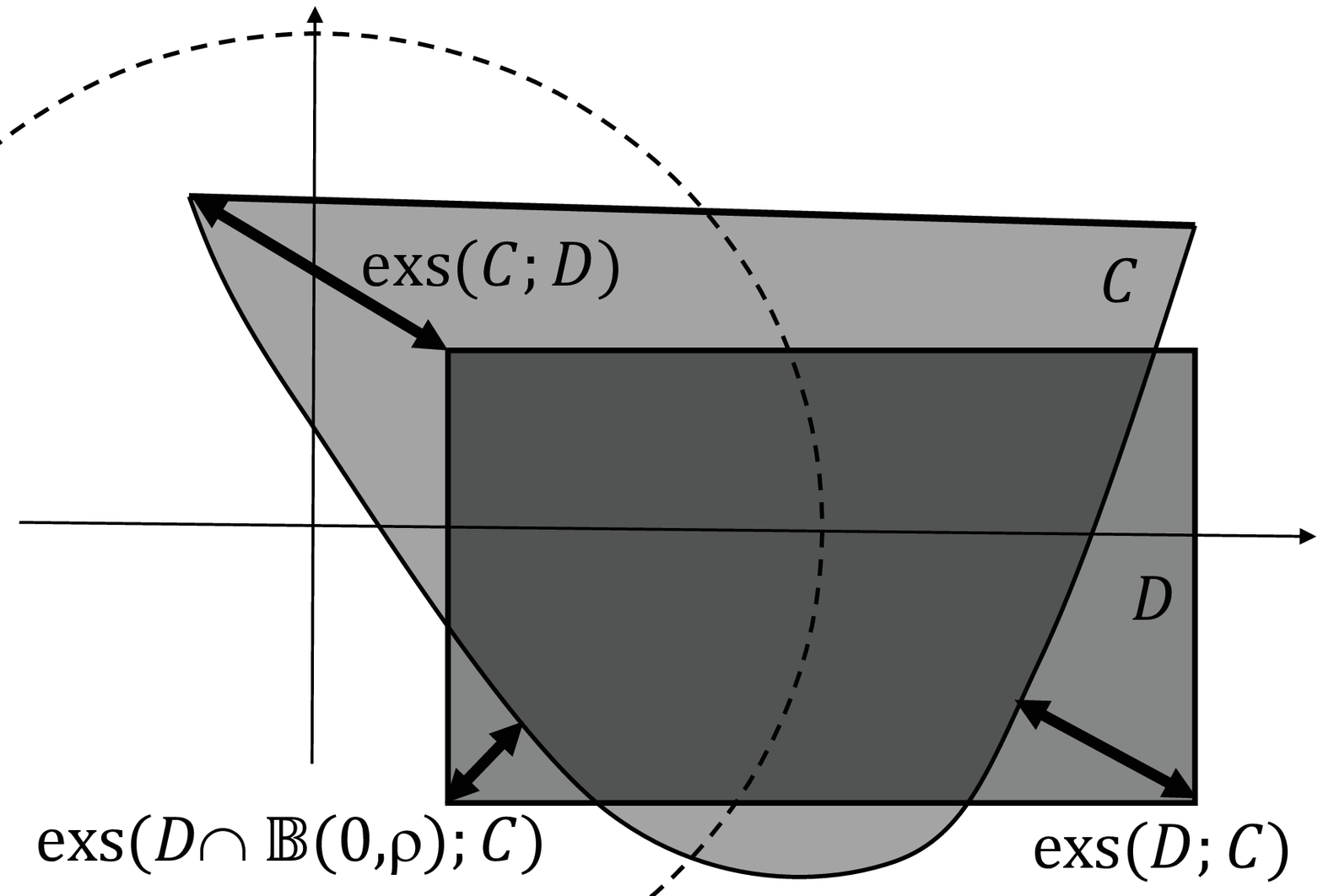}{2.6in} {The excess of $C$ over $D$ and vice versa. The truncated Hausdorff distance $\hatsetd_\rho(C,D) = \max\{\exs( C \cap \ball(0,\rho); D \big), ~\exs( D \cap \ball(0,\rho); C)\}$ considers only excess of the portion of a set inside a ball with radius $\rho$ (dashed circle). Here, the Euclidean norm is assumed, but any other can be employed as well.} {fig_excess}

For any $\rho\geq 0$ and $C,D\subset\reals^n$, the {\it truncated Hausdorff distance} between $C$ and $D$ is
\[
\hatsetd_\rho(C,D) = \max\Big\{\exs\big( C \cap \ball(0,\rho); D \big), ~~\exs\big( D \cap \ball(0,\rho); C\big)\Big\}.
\]
Figure \ref{fig_excess} exemplifies the effect of truncation:  $\exs( C \cap \ball(0,\rho); D )=\exs( C ; D )$ in this case because the portion of $C$ furthest away from $D$ is contained in $\ball(0,\rho)$, while $\exs( D \cap \ball(0,\rho); C )$ is quite a bit smaller than $\exs( D ; C )$. Still, $\hatsetd_\rho(C,D) = \exs(C;D)$ in this figure. Generally, $\hatsetd_\rho(C,D)$ is finite for any nonempty $C,D$ and $\rho\in \reals_+$.

The truncated Hausdorff distance characterizes set-convergence: Given $\bar\rho\geq 0$ and closed $C\subset\reals^n$,
\begin{equation}\label{eqn:setconCharact}
C^\nu\sto C ~\Longleftrightarrow~ \hatsetd_\rho(C^\nu, C)\to 0 \mbox{ for all } \rho\in [\bar\rho, \infty).
\end{equation}
The need for considering all $\rho$ sufficiently large is caused by the fact that the truncated Hausdorff distance would ``miss'' differences between $C^\nu$ and $C$ if they only occur further away from the origin than $\rho$ (in the norm used); see \cite[Theorem 4.36]{VaAn} for a proof and details\footnote{Any norm can be used to underpin the right-hand side of \eqref{eqn:setconCharact}, but this is also the case for the left-hand side. Even though Definition \ref{dSetLimit} is stated in terms of $\|\cdot\|_2$, any norm could have been adopted because, by Proposition \ref{pSetConvergence}, it's only the topology on $\reals^n$ that matters and all norms produce the same topology in finite dimensions.}.

In particular, we are interested in distances between epigraphs as these relate to distances between minimizers and minimum values. For $f,g:\reals^n\to \Reals$, $\epi f$ and $\epi g$ are subsets of $\reals^{n+1}$ and hence require a norm on $\reals^{n+1}$. For any norm $\|\cdot\|$ on $\reals^n$, a main choice for norm on $\reals^{n+1}$ is
\begin{equation}\label{eqn:prodnorm}
\max\big\{\|x-\bar x\|, |\alpha-\bar \alpha|\big\} \mbox{ for } (x,\alpha), (\bar x, \bar \alpha)\in \reals^n\times \reals,
\end{equation}
which typically results in the simplest formulae. In the context of minimization, we then obtain the following bounds; see \cite[Proposition 2.2]{Royset.20b}.

\begin{theorem}\label{thm:approxoptimalvalue}{\rm (approximation of minima and near-minimizers).}
  For $f,g:\reals^n\to\Reals$, $\epsilon,\rho\in \reals_+$, and any norm on $\reals^n$ with \eqref{eqn:prodnorm} for $\reals^{n+1}$, suppose that
  \begin{equation*}
    \ninf f, \ninf g \in [-\rho,\rho-\epsilon), ~~\nargmin f \cap \ball(0,\rho)\neq \emptyset, ~~\nargmin g \cap \ball(0,\rho)\neq \emptyset.
  \end{equation*}
  Then,
  \begin{align*}
  \big|\ninf f - \ninf g\big| &\leq \hatsetd_\rho(\epi f, \epi g)\\
  \exs\big(\epsilon\mbox{-}\nargmin g \cap \ball(0,\rho); ~\delta\mbox{-}\nargmin f\big) & \leq \hatsetd_\rho(\epi f, \epi g) \mbox{ for } \delta > \epsilon + 2\hatsetd_\rho(\epi f, \epi g).
  \end{align*}
\end{theorem}
The theorem holds for nearly arbitrary functions as long as minimizers are attained with finite minimum values and $\rho$ is large enough. The bounds are sharp as discussed in \cite{Royset.18,Royset.20b}.

Computation of the truncated Hausdorff distance is supported by the following condition that originated with Kenmochi \cite{Kenmochi.74} and fully realized by Attouch and Wets \cite{AttouchWets.91}. For the present form, see \cite[Proposition 4.1]{Royset.20b}. We denote {\it lower level-sets} of $f:\reals^n\to \Reals$ by $\{f\leq \delta\} = \{x\in \reals^n~|~f(x)\leq \delta\}$ for any $\delta \in \Reals$.

\begin{proposition}\label{prop:altform}{\rm (Kenmochi condition).}
For $f,g:\reals^n\to\Reals$, both with nonempty epigraphs, $\rho\in \reals_+$, and any norm on $\reals^n$ with \eqref{eqn:prodnorm} for $\reals^{n+1}$, we have
\begin{align*}
\hatsetd_\rho(\epi f, \epi g) =  &\ninf\Big\{ \eta\geq 0 ~\Big|~\inf_{x\in \ball(\bar x,\eta)} g(x) \leq \max\{f(\bar x), -\rho\} + \eta, ~~\forall \bar x\in \{f\leq \rho\} \cap \ball(0,\rho)\\
                            &~~~~~~~~~~~~~~~~~\inf_{x\in \ball(\bar x,\eta)} f(x) \leq \max\{g(\bar x), -\rho\} + \eta, ~~\forall \bar x\in \{g\leq \rho\} \cap \ball(0,\rho) \Big\}.
\end{align*}
\end{proposition}

From this fact, we can develop several useful bounds on the truncated Hausdorff distance. We say that a function $f:\reals^n\to\reals$ is {\it Lipschitz continuous} with modulus $\kappa:\reals_+\to \reals_+$ if
  \begin{equation*}
    \big|f(x) - f(\bar x)\big| \leq \kappa(\rho)\|x-\bar x\|_2 \mbox{ for all } \rho\in \reals_+ \mbox{ and } x,\bar x\in \ball_2(0,\rho) = \big\{x\in\reals^n~\big|~\|x\|_2\leq \rho\big\}.
  \end{equation*}
Note that the modulus is permitted to depend on the size of the ball on which the condition holds.

\begin{proposition}\label{prop:lipschitz} {\rm (distance between composite functions).}
   For proper $f_0,g_0:\reals^n\to \Reals$, consider $f = f_0+h\circ F$ and $g=g_0+h\circ G$, where $F = (f_1, \dots, f_m)$, $G=(g_1, \dots, g_m)$,
   \begin{enumerate}[{\rm (a)}]

    \item $h:\reals^m\to \reals$ is Lipschitz continuous with modulus $\kappa:\reals_+\to \reals_+$, and

    \item $f_i,g_i:\reals^n\to \reals$, $i=1, \dots, m$, are Lipschitz continuous with common modulus $\lambda:\reals_+\to \reals_+$.

   \end{enumerate}
   If $\|\cdot\|_2$ is the norm on $\reals^n$ and \eqref{eqn:prodnorm} with $\|\cdot\| = \|\cdot\|_2$ is the norm on $\reals^{n+1}$, then for any $\rho\in \reals_+$,
  \[
     \hatsetd_\rho(\epi f, \epi g) \leq \big(1+\sqrt{m}\kappa(\rho^*)\lambda(\hat \rho)\big)\hatsetd_{\bar\rho}(f_0, g_0) + \kappa(\rho^*)\nsup_{x\in \ball_2(0,\rho)} \big\|F(x) - G(x)\big\|_2,
  \]
where $\bar\rho \geq \rho + \max\{\nsup_{x\in\ball_2(0,\rho)} |h(F(x))|, \nsup_{x\in\ball_2(0,\rho)} |h(G(x))|\}$, $\hat\rho > \rho + \hatsetd_{\bar\rho} (\epi f_0, \epi g_0)$, and $\rho^* \geq \max\{\sup_{x\in \ball_2(0,\hat\rho)} \|F(x)\|_2,  \sup_{x\in \ball_2(0,\hat\rho)} \|G(x)\|_2 \}$.
\end{proposition}
\state Proof. Our goal is to establish the two conditions on the right-hand side in Proposition \ref{prop:altform}. Let $\mathbb{S}(\bar\rho) = \{(x,\alpha)\in\reals^n\times\reals~|~\|x\|_2\leq \bar\rho, |\alpha|\leq\bar\rho\}$ be a centered ball using the assumed norm on $\reals^{n+1}$. Set $\eta = \hatsetd_{\bar\rho}(\epi f_0, \epi g_0)$, $\epsilon \in (0, \hat\rho -\rho - \eta)$, and $\bar x\in \{f \leq \rho\} \cap \ball_2(0,\rho)$. Then, $f_0(\bar x) \leq \rho -h(F(\bar x)) \leq \bar\rho$.

First, suppose that $f_0(\bar x) \geq -\bar \rho$ so that $(\bar x,f_0(\bar x)) \in \epi f_0 \cap \mathbb{S}(\bar \rho)$. Consequently, there exists $(\hat x, \hat \alpha) \in \epi g_0$ with $\|\bar x -\hat x\|_2 \leq \eta+\epsilon$ and $|\hat\alpha - f_0(\bar x)|\leq \eta+ \epsilon$. Thus, $g_0(\hat x) \leq \hat \alpha \leq f_0(\bar x) + \eta+\epsilon$ and
\begin{align*}
 & \inf_{x\in \ball(\bar x,\eta+\epsilon)} g(x) \leq g_0(\hat x) + h\big(G(\hat x)\big) = g_0(\hat x) + h\big(F(\bar x)\big) + h\big(G(\hat x)\big) - h\big(G(\bar x)\big) + h\big(G(\bar x)\big) - h\big(F(\bar x)\big)\\
  & \leq f_0(\bar x)   + \eta + \epsilon + h\big(F(\bar x)\big) + \kappa(\rho^*)\big\|G(\hat x) - G(\bar x)\big\|_2 + \kappa(\rho^*)\big\|F(\bar x) - G(\bar x)\big\|_2\\
  & \leq f_0(\bar x) + h\big(F(\bar x)\big) + \eta + \epsilon+ \kappa(\rho^*)\sqrt{m}\max_{i=1, \dots, m} \big|g_i(\hat x) - g_i(\bar x)\big| + \kappa(\rho^*)\sup_{x\in \ball_2(0,\rho)}\big\|F(x) - G(x)\big\|_2\\
  & \leq f_0(\bar x) + h\big(F(\bar x)\big) + \eta + \epsilon+ \kappa(\rho^*)\sqrt{m}\lambda(\hat\rho)(\eta+\epsilon) + \kappa(\rho^*)\sup_{x\in \ball_2(0,\rho)}\big\|F(x) - G(x)\big\|_2\\
  & \leq \max\big\{f_0(\bar x) + h\big(F(\bar x)\big), -\rho\big\} + \big(1 + \kappa(\rho^*)\sqrt{m}\lambda(\hat\rho)\big)(\eta+\epsilon) + \kappa(\rho^*)\sup_{x\in \ball_2(0,\rho)}\big\|F(x) - G(x)\big\|_2.
\end{align*}

Second, suppose that $f_0(\bar x) < -\bar \rho$. Then, $(\bar x,-\bar \rho) \in \epi f_0 \cap \mathbb{S}(\bar \rho)$ and there exists $(\hat x, \hat \alpha) \in \epi g_0$ with $\|\bar x -\hat x\|_2 \leq \eta+\epsilon$ and $|\hat\alpha +\bar \rho|\leq \eta+ \epsilon$. Thus, $g_0(\hat x) \leq \hat \alpha \leq -\bar\rho+\eta+\epsilon$ and, similar to above,
\begin{align*}
  &\inf_{x\in \ball(\bar x,\eta+\epsilon)} g(x)  \leq -\bar\rho + h\big(F(\bar x)\big) + \eta + \epsilon+ \kappa(\rho^*)\sqrt{m}\lambda(\hat\rho)(\eta+\epsilon) + \kappa(\rho^*)\sup_{x\in \ball_2(0,\rho)}\big\|F(x) - G(x)\big\|_2\\
  &~~~~~~~~~~~ \leq \max\big\{f_0(\bar x) + h\big(F(\bar x)\big), -\rho\big\} + \big(1 + \kappa(\rho^*)\sqrt{m}\lambda(\hat\rho)\big)(\eta+\epsilon) + \kappa(\rho^*)\sup_{x\in \ball_2(0,\rho)}\big\|F(x) - G(x)\big\|_2.
\end{align*}
The second inequality follows because $- \bar\rho + h(F(\bar x)) \leq - \bar\rho + \nsup_{x\in\ball_2(0,\rho)}|h(F(x))| \leq -\rho$. Thus, in both cases, we obtain the same upper bound on $\inf_{x\in \ball(\bar x,\eta+\epsilon)} g(x)$. Repeating these arguments with the roles of $f$ and $g$ switched, we obtain the result via Proposition \ref{prop:altform} after letting $\epsilon$ tend to zero.\eop

For other rules to compute the truncated Hausdorff distance, we refer to \cite{AzePenot.90b,AttouchWets.91,VaAn,Royset.18,Royset.20b}.

We next turn to solution errors for generalized equations. If $S:\reals^n\tto\reals^m$, then $\gph S \subset \reals^n\times \reals^m$ and the truncated Hausdorff distance between such sets requires a norm on $\reals^n\times \reals^m$. Given any norm $\|\cdot\|_a$ on $\reals^n$ and any norm $\|\cdot\|_b$ on $\reals^m$,
a main choice is
\begin{equation}\label{eqn:prodnormsetval}
\max\big\{ \|x\|_a, \|y\|_b \big\} \mbox{ for } (x,y)\in \reals^n\times\reals^m.
\end{equation}
This leads to the following bound on near-solutions proven in \cite[Theorem 5.1]{Royset.20b}.

\begin{theorem}\label{thm:approxgeneralizedequations}{\rm (near-solutions of generalized equations).} For $S,T:\reals^n\tto \reals^m$, with nonempty graphs, and any norms on $\reals^n$ and $\reals^m$ with \eqref{eqn:prodnormsetval} for $\reals^n\times\reals^m$, suppose that $0 \leq\epsilon\leq \rho< \infty$ and $\bar y\in \ball(0,\rho-\epsilon)$. If $\delta > \epsilon + \hatsetd_\rho(\gph S, \gph T)$, then
\[
\exs\Big(S^{-1}\big(\ball(\bar y, \epsilon)\big)\cap \ball(0,\rho); ~T^{-1}\big(\ball(\bar y, \delta)\big)    \Big) \leq \hatsetd_\rho(\gph S, \gph T).
\]
If $\gph T$ is closed, then the result also holds with $\delta = \epsilon + \hatsetd_\rho(\gph S, \gph T)$.
\end{theorem}

The theorem is sharp. For example, consider $S,T:\reals\tto\reals$ with $S(x) = [x,\infty)$ when $x\in [0,1]$ and $S(x) = \emptyset$ otherwise; and $T(x) = (1,\infty)$ when $x\in [1,2]$ and $T(x) = \emptyset$ otherwise. Then for $\rho \geq 0$, $\hatsetd_\rho(\gph S, \gph T) = 1$, $S^{-1}(0) = \{0\}$, $T^{-1}(\delta)=[1,2]$, and $\exs(S^{-1}(0)\cap \ball(0,\rho); T^{-1}(\ball(0,\delta)) = 1$  when $\delta > 1$. When $\delta \leq 1$, the excess becomes infinity because $T^{-1}(\delta) = \emptyset$. If $T$ is modified to having $T(x) = [1,\infty)$ for $x\in [1,2]$, then $\delta = 1$ gives an excess of one.

\begin{example}{\rm (stationarity of composite functions).} For a proper lsc function $\phi:\reals^m\to \Reals$ and a smooth mapping $F:\reals^n\to \reals^m$, the composite function $\phi\circ F$ has $0 \in \nabla F(x)^\top \partial \phi(F(x))$ as a necessary optimality condition under standard assumptions \cite[Theorem 10.6]{VaAn}, where the $m\times n$-matrix $\nabla F(x)$ is the Jacobian of $F$ at $x$. Instead of solving this generalized equation, we may have to settle for an approximation. The error caused by this switch can be quantified using Theorem \ref{thm:approxgeneralizedequations}.
\end{example}
\state Detail. By introducing auxiliary vectors $y,z\in\reals^m$, the optimality condition is equivalently stated in terms of the set-valued mapping $S:\reals^n\times\reals^m\times\reals^m\tto \reals^m\times\reals^m\times\reals^n$ as $0 \in S(x,y,z)$, with
\begin{equation*}\label{eqn:compositionMapping}
  S(x,y,z) = \big\{F(x) - z\big\} \times  \big\{\partial\phi(z) - y\big\} \times \big\{\nabla F(x)^\top y\big\}.
\end{equation*}
Since $0\in S(x,y,z)$ is also an optimality condition for the problem of minimizing $\phi(z)$ subject to $F(x) = z$, we see that $y$ can be interpreted as a multiplier vector and $z$ as representing feasibility.

An approximating function $\psi\circ G$ expressed in terms of proper lsc $\psi:\reals^m\to \Reals$ and smooth $G:\reals^n\to\reals^m$ has parallel optimality conditions: $0\in T(x,y,z)$, where
\begin{equation*}\label{eqn:compositionMapping2}
  T(x,y,z) = \big\{G(x) - z\big\} \times \big\{\partial \psi(z) - y\big\} \times \big\{\nabla G(x)^\top y\big\}.
\end{equation*}
In view of Theorem \ref{thm:approxgeneralizedequations}, a bound on $\hatsetd_\rho(\gph S, \gph T)$ leads to an estimate of the change in near-stationary points as we pass from $\phi\circ F$ to $\psi\circ G$. 

The ``input space'' for $S,T$ is $\reals^n\times\reals^m\times\reals^m$ on which we adopt the norm $\max\{\|x\|_2,\|y\|_2,\|z\|_2\}$ for $(x,y,z) \in \reals^n\times\reals^m\times\reals^m$. The ``output space'' $\reals^m\times\reals^m\times\reals^n$ is assigned the norm $\max\{\|u\|_2,\|v\|_2,\|w\|_2\}$ for $(u,v,w) \in \reals^m\times\reals^m\times\reals^n$. Then, by \cite[Theorem 5.3]{Royset.20b}, for $\rho\in \reals_+$,
\[
\hatsetd_\rho(\gph S, \gph T) \leq \sup_{\|x\|_2\leq\rho} \max\Big\{ \big\|G(x) - F(x)\big\|_2 + \hatsetd_{2\rho} \big(\gph \partial \phi, \gph \partial \psi\big), ~\rho\big\|\nabla G(x) - \nabla F(x)\big\|_F \Big\},
\]
where $\|\cdot\|_F$ is the Frobenius norm. In calculating distances between $\gph \partial \phi$ and $\gph \partial \psi$, we adopt the norm $\max\{\|z\|_2,\|y\|_2\}$ for $(z,y)\in \reals^m\times \reals^m$.\eop

\section{Extensions}\label{sec:concl}

The quantification of set-convergence hints to the possibility of constructing a metric on spaces of sets and thereby enabling general metric space theory. Although $\hatsetd_\rho$ isn't a metric because it fails the triangle inequality and there is also the issue with discrepancies beyond $\ball(0,\rho)$, the closely related Attouch-Wets distance\footnote{This distance emerged from the work of Attouch and Wets in the early 1990s; see \cite[Section 4.I]{VaAn} where it's called the integrated set distance.} metrizes set-convergence on the space of nonempty closed sets. The resulting metric space is complete, separable, and finitely compact. The latter means that every closed ball in that space is compact. This points to the following fact available directly from the properties of set-convergence: For every sequence $\{C^\nu\subset\reals^n, n\in\nats\}$, either $C^\nu \sto \emptyset$ or there exists $N\in \cN_\infty^\grill$ and a nonempty set $C$ such that $\{C^\nu, \nu\in N\} \sto C$. This result can be traced back to Zoretti in 1909, a student of Painlev\'{e}, but here presented in the form found in \cite[Theorem 4.18]{VaAn}.

When restricted to epigraphs of lsc functions with nonempty domains, the Attouch-Wets distance furnishes a metric on the space of such functions with the property that convergence of functions implies convergence of the corresponding minimizers and minimum values in the sense of Theorem \ref{tEpiConvConsequences}(a,b,c). This can be utilized to analyze infinite-dimensional problems in nonparametric statistics and elsewhere; see \cite{RoysetWets.15b,RoysetWets.19,Royset.20}. The compactness of closed and bounded subsets of such spaces trivializes questions of existence of solutions, results in finite covering numbers \cite{Royset.20},  and facilitates the construction of approximating functions that are computationally attractive \cite{RoysetWets.15b,Royset.18,Royset.20}.

The extension of set-convergence and its quantification from subsets of $\reals^n$ to those of general metric spaces is trivial as the definitions only require a notion of distance between points. The truncated Hausdorff distance is then defined relative to an arbitrary center point; for $\reals^n$ we simply adopt the origin, which gives rise to $\ball(0,\rho)$ in the definition. In fact, even on $\reals^n$, it might be beneficial to move the center point elsewhere to better capture the scale of a problem; see \cite{Royset.18,Royset.20}. A complication beyond $\reals^n$ is that \eqref{eqn:setconCharact} may hold only with $\Longleftarrow$; see \cite{Royset.18} for an introduction and \cite{Beer.93} for details. Nevertheless, set-convergence remains an important tool in the infinite-dimensional setting; see for example \cite{AttouchWets.93a,Polak.97,KorfWets.01,AttouchButtazzoMichaille.14,PhelpsRoysetGong.15,Royset.18}.\\

\noindent {\bf Acknowledgement.} This work is supported by ONR (Operations Research) under  N0001420WX00519 and AFOSR (Optimization and Discrete Mathematics) under  F4FGA08272G001.

\bibliographystyle{plain}
\bibliography{refs}

\end{document}

%% file: mac.tex
\usepackage{theorem}
\usepackage{enumerate}
\usepackage{array}
\usepackage[reqno]{amsmath}
\usepackage{amssymb}
\usepackage{mathtools}
\usepackage{latexsym}
\usepackage{makeidx}
\usepackage{fancybox}
\input epsf.sty
\usepackage[dvips]{graphicx,color}
\usepackage{showlabels}
\usepackage{subcaption}

\newcommand{\drawing}[4]{
	\begin{figure}[!hbt]
	\begin{center}
	\leavevmode
	\epsfxsize=#2
	\epsfbox{#1}
	\caption{\small #3}
	\label{#4}
	\end{center}
	\end{figure}}

\theoremstyle{change}
\newtheorem{proclaim}{PROCLAIM}[section]
\newtheorem{theorem}[proclaim]{Theorem}
\newtheorem{definition}[proclaim]{Definition}

\newtheorem{proposition}[proclaim]{Proposition}

\newtheorem{example}[proclaim]{Example}

\numberwithin{equation}{section}

\outer\def\proclaim #1. #2\par{\medbreak \noindent{\bf#1.\enspace}{\sl#2}\par
  \ifdim\lastskip<\medskipamount
  \removelastskip\penalty55\medskip\fi}
\def\state #1. { \noindent{\bf#1.\enspace}}
\def\algo #1. { \noindent{\bf#1.\enspace}}

\DeclareMathOperator{\con}{con}

\DeclareMathOperator{\dist}{dist}

\DeclareMathOperator{\exs}{exs}
\DeclareMathOperator{\dom}{dom}
\DeclareMathOperator{\gph}{gph}
\DeclareMathOperator{\epi}{epi}
\DeclareMathOperator{\hypo}{hypo}
\DeclareMathOperator{\nt}{int}

\DeclareMathOperator{\prj}{prj}

\newcommand{\comp}{\,{\raise 1pt \hbox{$\scriptstyle\circ$}}\,}

\newcommand{\reals}{\mathbb{R}}
\newcommand{\Reals}{\overline{\mathbb{R}}}
\newcommand{\natnums}{{{\rm l} \kern -.13em {\rm N} }}
\newcommand{\nats}{\mathbb{N}}
\newcommand{\snats}{{I\kern -.29em N}}
\newcommand{\rats}{{Q\kern -.64em \raise 1pt \hbox{$\scriptstyle |$}\;\,}}
\newcommand{\srats}
	{{Q\kern -.56em \raise 1.2pt \hbox{$\scriptscriptstyle /$}\,}}
\newcommand{\ints}{Z\kern -.46em Z}
\newcommand{\ball}{\mathbb{B}}
\newcommand{\pluss}{\hskip1pt \raise1pt\vbox{\hrule width6pt \vskip1pt \hrule
                    width6pt} \kern-4pt{\lower1pt\hbox{\vrule height6pt
		    \kern1pt\vrule height6pt}}\hskip5pt}
\newcommand{\eop}
	{\hfill{$\vcenter{\hrule height1pt \hbox{\vrule width1pt height5pt
   	 \kern5pt \vrule width1pt} \hrule height1pt}$} \medskip}
\newcommand{\half}
	{{\raisebox{1pt}{$\frac{1}{2}$}}}

\newcommand{\setd}{{ d \kern -.15em l}}
\newcommand{\hatsetd}{ d \hat{\kern -.15em l }}

\renewcommand{\epsilon}{\varepsilon}
\renewcommand{\phi}{\varphi}

\newcommand{\tto}{\;{\lower 1pt \hbox{$\rightarrow$}}\kern -12pt
           \hbox{\raise 2.5pt \hbox{$\rightarrow$}}\;}
\newcommand{\overto}[1]{\,{\raise 0pt\hbox{$\rightarrow$}}\kern -9pt
     \hbox{\lower 3pt \hbox{$\scriptscriptstyle#1$}}\hskip6pt}
\newcommand{\underto}[1]{\,{\lower 1pt\hbox{$\rightarrow$}}\kern -9pt
     \hbox{\raise 4pt \hbox{$\,\scriptscriptstyle#1$}}\hskip7pt}
\newcommand{\bigoverto}[1]{{\raise 0pt\hbox{$\,\longrightarrow$}}\kern -16pt
     \hbox{\lower 3pt \hbox{$\scriptscriptstyle#1$}}\hskip4pt}
\newcommand{\bigunderto}[1]{\,{\lower 1pt\hbox{$\longrightarrow$}}\kern -16pt
     \hbox{\raise 4pt \hbox{$\,\scriptscriptstyle#1$}}\hskip6pt}
\newcommand{\bigbigto}[2]{\,{\raise 0pt\hbox{$\,\longrightarrow$}}\kern -16pt
     \hbox{\lower 3pt \hbox{$\scriptscriptstyle#2$}}\kern -10pt
     \hbox{\raise 4pt \hbox{$\,\scriptscriptstyle#1$}}\hskip7pt}
\newcommand{\downto}{{\raise 1pt \hbox{$\scriptscriptstyle \,\searrow\,$}}}
\newcommand{\upto}{{\raise 1pt \hbox{$\scriptscriptstyle \,\nearrow\,$}}}

\newcommand{\notimply}
	{\quad\hbox{$\Longrightarrow \kern -14pt {/}$}\hskip6pt\quad}

\newcommand{\lto}{\,{\lower 1pt\hbox{$\rightarrow$}}\kern -10pt
     \hbox{\raise 4pt \hbox{$\, \scriptstyle l$}}\hskip7pt}
\newcommand{\eto}{\,{\lower 1pt\hbox{$\rightarrow$}}\kern -11pt
     \hbox{\raise 4pt \hbox{$\, \scriptstyle e$}}\hskip7pt}
\newcommand{\hto}{\,{\lower 1pt\hbox{$\rightarrow$}}\kern -11pt
     \hbox{\raise 4pt \hbox{$\, \scriptstyle h$}}\hskip7pt}
\newcommand{\pto}{\,{\lower 1pt\hbox{$\rightarrow$}}\kern -11pt
     \hbox{\raise 4.5pt \hbox{$\, \scriptstyle p$}}\hskip7pt}
\newcommand{\cto}{\,{\lower 1pt\hbox{$\rightarrow$}}\kern -11pt
     \hbox{\raise 4pt \hbox{$\, \scriptstyle c$}}\hskip7pt}
\newcommand{\gto}{\,{\lower 1pt\hbox{$\rightarrow$}}\kern -11pt
     \hbox{\raise 4.5pt \hbox{$\, \scriptstyle g$}}\hskip7pt}
\newcommand{\sto}{\,{\lower 1pt\hbox{$\rightarrow$}}\kern -11pt
     \hbox{\raise 4pt \hbox{$\, \scriptstyle s$}}\hskip7pt}
\newcommand{\awto}{\,{\lower 1pt\hbox{$\rightarrow$}}\kern -15pt
     \hbox{\raise 4pt \hbox{$\, \scriptstyle aw$}}\hskip7pt}
\def\Nto{\,{\raise 1pt\hbox{$\rightarrow$}}\kern -13pt
     \hbox{\lower 3pt \hbox{$\, \scriptstyle N$}}\hskip7pt}
\def\Cto{\,{\raise 1pt\hbox{$\rightarrow$}}\kern -14pt
     \hbox{\lower 3pt \hbox{$\, \scriptstyle C$}}\hskip7pt}
\def\fto{\,{\raise 1pt\hbox{$\rightarrow$}}\kern -14pt
     \hbox{\lower 3pt \hbox{$\, \scriptstyle f$}}\hskip7pt}

\newcommand{\low}[1]{{\lower1pt \hbox{$\scriptstyle #1$}}}
\newcommand{\loww}[1]{{\lower2pt \hbox{$\scriptstyle #1$}}}
\newcommand{\high}[1]{{\raise1pt \hbox{$\scriptstyle #1$}}}

\newcommand{\cN}{{\cal N}}

\newcommand{\cP}{{\cal P}}

\newcommand{\nsum}{\mathop{\sum}\nolimits}

\newcommand{\nInnLim}{\mathop{\rm LimInn}\nolimits}
\newcommand{\nOutLim}{\mathop{\rm Lim\hspace{-0.01cm}Out}\nolimits}
\newcommand{\nLim}{\mathop{\rm Lim}\nolimits}

\newcommand{\nlim}{\mathop{\rm lim}\nolimits}
\newcommand{\nliminf}{\mathop{\rm liminf}\nolimits}

\newcommand{\nlimsup}{\mathop{\rm limsup}\nolimits}
\newcommand{\ninf}{\mathop{\rm inf}\nolimits}
\newcommand{\nsup}{\mathop{\rm sup}\nolimits}

\newcommand{\nnmin}{\mathop{\rm minimize}}

\newcommand{\nargmin}{\mathop{\rm argmin}\nolimits}

\newcommand{\lwdy}[2]{\mathrel{\mathop
        {\raisebox{0.1ex}{\null$#1$}}{\hbox{\kern -1.0em
	{\raisebox{-0.8ex}{$\scriptstyle{\;\to #2}$}}}}}}
\newcommand{\lwwdy}[2]{\mathrel{\mathop
        {\raisebox{0.2ex}{\null$#1$}}{\hbox{\kern -1.0em
	{\raisebox{-1.1ex}{$\scriptstyle{\;\to #2}$}}}}}}
\newcommand{\slwdy}[2]{\scriptsize{{\mathrel{\mathop
        {\raisebox{0.1ex}{\null$#1$}}{\hbox{\kern -1.0em
	{\raisebox{-0.8ex}{$\scriptstyle{\;\to #2}$}}}}}}}}
\newcommand{\slwwdy}[2]{\scriptsize{{\mathrel{\mathop
        {\raisebox{0.2ex}{\null$#1$}}{\hbox{\kern -1.0em
	{\raisebox{-1.1ex}{$\scriptstyle{\;\to #2}$}}}}}}}}

\definecolor{lightgray}{gray}{0.75}
\definecolor{myred}{rgb}{0.55,0,0}
\definecolor{myblue}{rgb}{0,0,0.5} 
\definecolor{mygreen}{rgb}{0,0.5,0} 
\definecolor{purple}{rgb}{0.5,0,0.5} 
\definecolor{turq}{rgb}{0,0.805,0.816} 
\definecolor{maroon}{rgb}{0.51,0,0}
\definecolor{MAROON}{rgb}{0.51,0,0}
\definecolor{redor}{rgb}{0.78,0.078,0.078}
\definecolor{dgreen}{rgb}{0,0.3,0}

\newcommand{\grill}{{\scriptscriptstyle\#}}
\newcommand{\bcdot}{\,{\raise .2ex \hbox{$\centerdot$}}\,}

